\documentclass[11pt]{amsart}

\usepackage{amsmath,amssymb,amsthm}
\usepackage{graphicx}
\usepackage{amssymb}

\numberwithin{equation}{section}
\newtheorem{thm}{Theorem}
\newtheorem{lemma}[thm]{Lemma}

\newtheorem{corollary}[thm]{Corollary}
\newtheorem{proposition}[thm]{Proposition}
\newtheorem{definition}{Definition}{\rm}
\newtheorem{example}{Example}{\rm}

{\rm}

\def\beq{\begin{equation} }
\def\eeq{\end{equation} }

\def\C{\mathbb{C}}

\def\N{\mathbb{N}}
\def\R{\mathbb{R}}

\def\cA{\mathcal A}
\def\bbN{\mathbb N}
\def\cB{\mathcal B}
\def\cM{\mathcal M}

\def\cN{\mathcal N}
\def\cD{\mathcal D}
\def\cK{\mathcal K}
\def\cM{\mathcal M}
\def\cQ{\mathcal Q}
\def\Hh{\mathfrak H}

\def\R{\mathbb{R}}

\newcommand{\ol}[1]{\overline{#1}}

\begin{document}
    \title[Multivariate Determinateness]{Multivariate Determinateness}
    \author[Putinar]{Mihai Putinar}
    \address{Department of Mathematics, University of California at Santa Barbara,
    Santa Barbara, CA 93106, USA}
    \email{mputinar@math.ucsb.edu}
    \urladdr{www.math.ucsb.edu/\~{}mputinar}
    \author[Schm\"udgen]{Konrad Schm\"udgen}
    \address{Mathematical Institute, University of Leipzig,
Postfach 10 09 20, D-04009 Leipzig, Germany}
    \email{schmuedgen@math.uni-leipzig.de}
    \urladdr{www.math.uni-leipzig.de/\~{}schmuedgen}
    \dedicatory{To Ciprian Foias, on the occasion of his seventy-fifth birthday}
    \date{\today}
    \thanks{The first author was partially supported by the National Science Foundation grant
    DMS-0701094}
    \begin{abstract} The uniqueness question of the multivariate
    moment problem is studied by different methods: Hilbert space operators, complex
    function theory, polynomial approximation, disintegration, integral geometry.
    Most of the known results in the multi-dimensional case are reviewed and
    reproved,
    and a number of new determinacy criteria are developed.
     \end{abstract}

    \maketitle

\section{Introduction}

A moment problem represents a quintessential inverse problem,
continuously reborn and rediscovered on the basis of novel,
surprising applied necessities. In the restricted sense we adopt
in this article, a moment problem revolves around reconstructing,
approximating or only estimating a positive Borel measure on a
closed subset of the Euclidean space from its power moments, or a
finite selection of them. Among all mathematical aspects of moment
problems, the uniqueness, or, following circulating terminologies,
the determinateness or determinacy question (that is, under which
conditions a sequence of moments corresponds to a single measure)
has attracted analysts for more than a century. With all efforts,
this very topics remains largely unfinished in the case of
multivariate moment problems. Our paper aims at offering a
complete picture of the present status of the determinateness
question for multivariate moment problems.

Let $\mu$ be a positive measure on the real line, having all power
moments finite. The asymptotic expansion in wedges $0< \delta
<\arg z < \pi -\delta$ of the associated Cauchy transform
$$\int_\R\frac{d\mu(x)}{x-z}\sim -\frac{s_0}z-\frac{s_1}{z^2}
-\cdots,\quad\Im(z)>0,$$ and the well defined continued fraction
$$-\frac{s_0}z-\frac{s_1}{z^2}-\cdots=-\cfrac{s_0}{z-\alpha_0
-\cfrac{\beta_0^2}{z-\alpha_1-\cfrac{\beta_1^2}{z-\alpha_2-\
\cfrac{\beta_2^2}{\ddots}}}},\quad\alpha_k,\>\beta_k\in\R,$$
depend solely on the sequence of moments
$$s_k=\int_\R x^k\>d\mu(x),\quad k\ge0.$$
The moment sequence $\{ s_k; \ k \geq 0\}$ is called
\emph{determinate} if there exists a unique representing positive
measure.

Determinateness criteria for the above (so called Hamburger)
moment problem are classical. One can start with the observation
that two measures $\mu$ and $\nu$ with the same moments coincide
if and only if their Cauchy transforms are equal on a set of
uniqueness for analytic functions in the upper half-plane. For
instance the above moment problem is determinate precisely when
the associated continued fraction converges at non-real points
$z$. Originally, all determinateness criteria were obtained via
elaborated computations involving continued fractions; for an
excellent presentation of these aspects see the old monograph by
Perron \cite{Pe}, Hamburger's work \cite{Ham1,Ham2} or even the
original memoir by Stieltjes \cite{Stieltjes} which defeats time
by its universality and freshness.

The orthogonal polynomials of the first and second kind $P_k(z)$,
respectively $Q_k(z), \ k \geq 0,$ again depending only on the
moment sequence $(s_k)$, introduce into the picture the Jacobi
matrix $J$, a three diagonal infinite, formally symmetric matrix
representing the multiplication by $z$ in the basis $(P_k(z))$:
$$ zP_k(z) = b_{k-1}P_{k-1}(z) + c_k P_k(z) + b_k P_{k+1}(z),\ \
k\in \N.$$ The (unbounded) self-adjoint extensions $A$ of $J$ on a
possibly larger Hilbert space than $l^2(\N)$ give a
parametrization of all solutions to the moment problem, via the
resolvent identity
$$ \int_\R\frac{d\nu(x)}{x-z} =  \langle (A-z)^{-1} 1, 1
\rangle,\ \ \Im z>0.$$ From here one derives M. Riesz' criterion
\cite{Ri}: the moment problem associated to the measure $\mu$ is
determinate if and only if the constant function $1$ is in the
closure of the polynomial ideal $(x+z)$ in $L^2(\mu)$, for a
single non-real complex number $z$, and hence for all $z$. For
these aspects of the theory of moments the reader can consult the
monographs \cite{A,ST,Stone} and the recent survey \cite{Simon}.

The link between moment problems and complex analytic function
theory is not less exciting and rewarding. It was Carleman
\cite{Carleman} who approached Stieljes moment problem (on the
semi-axis) with Laplace transform methods. Specifically, if $\mu$
is a positive measure on $\R_+ = [0,\infty),$ then the analytic
continuation of the Laplace transform
$$ \tilde{\mu}(z) = \int_0^\infty e^{-xz} d\mu(x), \ \ \Re z >0,$$
provides another effective uniqueness criterion: the quasi-analyticity of
$\tilde{\mu}(x)-\tilde{\nu}(x), x\geq 0,$ where $\nu$ is another
measure possessing the same moments, yields the determinacy of the
original moment problem.

In another direction, we owe to Nevanlinna a thorough analysis of
indeterminate moment problems on the line, via a canonical
representation of the Cauchy transforms of all solutions $\nu$:
$$ \int_\R\frac{d\nu(x)}{x-z}  = - \frac{C(z) \Phi(z)+ A(z)}{D(z)\Phi(z)+B(z)},\ \ \Im z >0.$$
The entire functions $A, B, C, D$ depend only on the orthogoonal
polynomials, and hence only on the sequence of moments. The
parameter $\Phi$ runs over all analytic functions in the upper
half-plane with values in the closure of it in the Riemann sphere.
A canonical parametrization of all solutions $\nu$ of the moment
problem is thus derived, see for details \cite{A,Simon}.

In sharp contrast to the one-dimensional case the uniqueness
question for {\it multivariate} moment problems is much less
understood. Riesz' density criterion has several counterparts (see
e.g. \cite{Fug}); however, all of them provide only sufficient
uniqueness conditions. Hilbert space methods lead to rather strong
determinateness results (cf.\ \cite{Devinatz, Nussbaum}), but
again they remain far from being also necessary. Integral
transforms give uniqueness results via quasi-analyticity and the
geometry of the support of the measure can also play a decisive
role. Disintegration methods and integral geometry techniques
naturally come into discussion and allow  possible dimension
reductions. Integral geometry leads even to a parametrization of
all solutions of the moment problem supported by a convex wedge.

The aim of this paper is to develop in a self-contained text the
above ideas. We expose (almost) all known facts and prove a number
of new results related to the multi-dimensional determinacy
problem. While we will always take the one-dimensional case for
granted, we shall give complete proofs of all main theorems referring
to the multi-dimensional case. These proofs are inserted into a
new, more natural, context; many of them depart from the original ones and
in this way we have corrected some existent unnoticed errors in the
literature.

The contents of this paper is the following. Section 2 deals with three
classical integrals: the Laplace, Fantappi\`e and respectively
Cauchy transforms of a positive measure in $\R^d$. Their
asymptotic expansions depend on the moments of the original
measure, and a discussion of the correspondence between moment
sequences (including the uniqueness question) and the function
spaces where their transforms belong is included.

Section 3 is devoted to the Hilbert space approach to  moment
problems. We begin with the Gelfand-Naimark-Segal construction and
the spectral theorem for strongly commuting self-adjoint
extensions of the multiplication operators and treat various
operator-theoretic determinacy notions in this context. We propose
a new concept, called {\it strict determinacy}, by requiring
determinacy together with the density of polynomials in
$L^2(\mu)$.

In Section 4 we derive a series of  determinacy criteria based  on
polynomial or rational approximations. They  include
generalizations and applications of Petersen' s theorem
\cite{Petersen} which are useful for verifying determinacy in many
situations, but also an approximation result obtained in
\cite{PS}.

Section 5 is centered around partial determinateness. If
determinacy of "sufficently many" $1$-subsequences of a positive
semi-definite multi-sequence is assumed, then the multi-sequence
is even a moment sequence and determinate. Pioneering results
belong to A. Devinatz \cite{Devinatz}, G.I. Eskin \cite{Eskin}, J.
Friedrich \cite{Fried} and others. In this section we develop some
general results of this kind and offer new proofs, again relying
essentially on self-adjointness theory.

In Section 6 the partial determinacy results of Section 5 are used to derive another classical result due to A.E. Nussbaum \cite{Nussbaum} on quasi-analyticity: a positive semi-definite multi-sequence is a determinate moment sequence if all marginal sequences satisfy  Carleman's quasi-analyticity condition.

Section 7 is an essay  on how a parametrization of all solutions
to the moment problem on an octant $\R_+^d$ in $\R^d$ should look
like. To this aim we combine a natural geometric push forward on a
pencil of lines with a characterization of the Fantappi\`e
transforms of positive measures on $\R_+^d$ due to Henkin and
Shananin \cite{HS}.

In Section 8 we apply disintegration
techniques to the determinacy problem. Roughly speaking, we study how determinateness assumptions on the
factors appearing in the disintegration formula
 $$\int_X~f(x) ~d\nu(x) = \int_T ~d\mu(t) ~\int_X
~f(x)~d\lambda_t(x)$$ imply the determinateness
of the measure $\nu$. Our main result is a new reduction theorem which shifts the determinacy question to lower  dimensions. A number of corollaries and applications show the usefulness of this result.

Finally, in Section 9 we reproduce the main result from \cite{PS}
which gives geometric criteria of determinateness. That is, quite
unexpectedly, there are closed and unbounded real algebraic
subsets of $\R^d$ which support only determined moment problems.
Again, a variety of geometric examples of low degree/low dimension
illustrate this phenomenon.

As mentioned above, we have tried to present the status of the
art of, and to develop new results referring to the multivariate determinacy problem.
We hope that this article will stir interest and stimulate further research towards a
deeper understanding of the multi-dimensional moment problem.\\

\tableofcontents

 {\small {\it Notations and definitions:}

Let $\R[x]{=}\R[x_1,\dots,x_d]$ and $\C[x]{=}\C[x_1,\dots,x_d]$ denote the  real and complex polynomial $\ast$-algebras, respectively. For  $\alpha{=}(\alpha_1,\dots,\alpha_d) \in \N_0^d$ and $x = (x_1,...,x_d) \in
\R^d$ we use the multi-index notation  $x^\alpha {=}
x_1^{\alpha_1}...x_d^{\alpha_d}$, where $x_j^0:=1$. Let $\R_+ := [0,\infty)$.

Let $\cM(\R^d)$ be the set of positive Borel measures $\mu$ on $\R^d$ such that $x^\alpha \in L^1(\mu)$ for all $\alpha \in \N_0^d$. Measures in $\cM(\R^d)$ are also called {\it rapidly decaying at
infinity} measures  or measures which have all moments.
For a closed subset $K$ of $\R^d$, $\cM(K)$ denotes the set of $\mu \in
\cM(\R^d)$ such that ${\rm supp}~\mu \subseteq K$. Whenever we use
the symbol $\cM(K)$ we mean that $K$ is a closed subset of some space
$\R^d$.

For a measure $\mu \in \cM(\R^d)$, the {\it $\alpha$-th moment} of $\mu$ is the number $ s_\alpha = \int x^\alpha d\mu(x)$ and  the {\it moment functional} $L_\mu$
is the linear functional on $\C[x]$  defined by
$L_\mu(p)=\int p(x) d\mu(x)$ for $p \in \C[x]$. Throughout we freely use the one-to-one correspondence between multi-sequences $\{s_\alpha;\alpha \in \N_0^d\}$ and linear functionals $L$ on $\C[x]$ given by $L(x^\alpha)=s_\alpha, \alpha \in \N_0$. The {\it moment problem} is the question when a given multi-sequence $s=\{s_\alpha;\alpha \in \N_0^d\}$ is the moment sequence of some measure $\mu \in \cM(\R^d)$ or equivalently when a given linear funtional $L$ is the moment functional of some $\mu \in \cM(\R^d)$. A necessary condition for this is  that the sequence is {\it positive semi-definite} (that is, $\sum s_{\alpha +\beta} \xi_\alpha \overline{\xi_\beta } \geq 0$ for any finite complex sequence $\{\xi_\alpha;\alpha \in \N_0^d\}$) or equivalently the correponding functional $L$ is {\it positive} (that is, $L(p\overline{p}) \geq 0$ for all $p\in \C[x]$), but this condition is not sufficient when $d \geq 2$.

We say that two measures $\mu, \nu \in \cM(\R^d)$ are  {\it
equivalent} and write $\mu \cong \nu$ if they have the same
moments or equivalently if $L_\mu =L_\nu$. Let $V_\mu$ denote the
set of all measures $\nu \in \cM(\R^d)$ such that $\mu \cong \nu$.
The set $V_\mu$ is convex and compact in the weak-* topology;

A measure $\mu \in \cM(K)$ (and likewise its moment sequence $s$ and its moment functioal $L_\mu$) is called {\it determinate on $K$} if any other measure $\nu \in
\cM(K)$ such that $\mu \cong \nu$ is equal to $\mu$. Let
$V_\mu(K)$ be the set of all positive measures belonging to $V_\mu$
and having support contained in $K$.

We say that $\mu$ and $L_\mu$ are {\it strongly determinate} if
the associated multiplication operators $X_1,\dots,X_d$ by the coordinates  are essentially
self-adjoint on the Hilbert space completion of polynomials. We
call $\mu$ and $L_\mu$ {\it strictly determinate} if $\mu$ and
$L_\mu$ is determinate and the polynomials $\C[x]$ are
dense in $L^2(\mu)$. Finally we say that $\mu$ and $L_\mu$ are called {\it ultradeterminate} if
the polynomials
$\C[x]$ are dense in $L^2((1+||x||^2)\mu)$.}\\

{\bf Acknowledgement.} This work was started in June 2007, during
a visit of the first author at the University of Leipzig.
    He wishes to express his gratitude to his hosts for hospitality and optimal working conditions.

\section{Review of integral transforms}

We focus on the retrieval of a positive measure $\mu$ from its
moments $\{s_\alpha;\alpha \in \N_0^d\}$. It is well known that
this is an ill posed problem, and that only under additional
assumptions the representing measure $\mu$ associated to a
prescribed moment sequence is unique.

Traditionally, the retrieval of $\mu$ is achieved by the exact or
approximative inversion of certain integral transforms of $\mu$.
These transforms contain the moment data in a closed form or an
asymptotic expansion form. We discuss below three such classical
examples.

\subsection{Laplace transform} Suppose that $\mu \in \cM(\R_+^d)$. Then the function
$$ \tilde{\mu}(z) = \int_{\R_+^d} e^{- t \cdot z} d\mu(t),$$
is analytic in the tube domain $\Re z \in {\rm int} \R_+^d.$
Remark that
$$ \frac{\partial^\alpha}{\partial z^\alpha} \tilde{\mu}(z) = (-1)^{|\alpha|}
\int_{\R_+^d} t^\alpha e^{- t \cdot z} d\mu(t),$$ whence
$$ \sup_{\Re z \in {\rm int} \R_+^d}|\frac{\partial^\alpha}{\partial z^\alpha} \tilde{\mu}(z)| \leq s_\alpha, \ \
\alpha \in \N_0^d.$$ It is well known that the Laplace transform
is injective, that is $\tilde{\mu} = 0$ implies $\mu=0$, see for
instance \cite{Ehrn}.

Regarded as a function of  real variables, $\tilde{\mu} \in
{\mathcal C}^\infty (\R_+^d)$ and we have
$$ (-1)^{|\alpha|}   [\frac{\partial^\alpha}{\partial x^\alpha} \tilde{\mu}](0) = s_\alpha, \ \ \alpha \in \N_0^d.$$

The function $\tilde{\mu}$ is {\it completely monotonic}, in the sense that
$$ (-1)^{|\alpha|}   [\frac{\partial^\alpha}{\partial x^\alpha} \tilde{\mu}](x) \geq 0, \ \ x \in \R_+^d, \ \
\alpha \in \N_0^d.$$

The following remarkable characterization of Laplace transforms
goes back to the works of S. Bernstein, V. Gilbert and S. Bochner.
For the proof we refer to \cite{HS}.

\begin{thm}  A function $F \in {\mathcal C}^\infty (\R_+^d)$ is the Laplace transform of a
measure from $\cM(\R_+^d)$, if and
only if $F$ is completely monotonic.
\end{thm}

Variants of this result for exponentially decaying measures, or
simply for all positive measures on $R_+^d$ are also well known,
see \cite{HS}.

Thus, for the uniqueness problem, we are naturally led to consider
two representing measures $\mu_1, \mu_2$ and take
$$ f(x) = \tilde{\mu_1}(x)-\tilde{\mu_2}(x), \ \  x \in \R_+^d.$$
Then
\begin{enumerate}
\item $f \in  {\mathcal C}^\infty (\R_+^d), $

\item $ \sup_{x \in \R_+^d} |f^{(\alpha)}(x)| \leq 2 s_\alpha, \ \ \alpha \in \N_0^d,$

\item$ f^{(\alpha)}(0) = 0, \ \ \alpha \in \N_0^d.$
\end{enumerate}

Uniqueness holds if $f$ is identically zero, that is this function
belongs to a quasi-analytic class. This was the original path
followed by Carleman in his application of quasi-analyticity to
the study of the uniqueness in the one variable moment problem
\cite{Carleman}. In the same spirit, we will discuss at the end of
this section how quasi-analyticity conditions can be applied to
find uniqueness criteria in the multi-variate moment problem on
the positive octant.

\subsection{Fantappi\`e transform} Closer to Stieljes original study \cite{Stieltjes}
of the
moment problem on the semi-axis is the {\it Fantappi\`e transform:}
$$ {\mathcal F}(\mu)(p_0,p) = \int_{\R_+^d} \frac{d\mu(x)}{p_0 + p \cdot x},\ \ p_0>0,\  p \in
\R_+^d.$$

A similar result to the characterization of Laplace transforms holds in this case.

\begin{thm}\label{Fantappie}
 A function $F \in {\mathcal C}^\infty({\rm int}\R \times \R_+^d)$ represents the Fantappi\`e
transform of a measure from $\cM(\R_+^d)$ if and only if $F$ is completely monotone
and homogeneous of degree $-1$.
\end{thm}

For a proof and variants, see again \cite{HS}. Note that the
moments $s_\alpha$ appear in the asymptotic expansion
$$  {\mathcal F}(\mu)(p_0,p) \sim \sum_{\alpha \in \N^d}
\frac{(-1)^{|\alpha|}}{p_0^{|\alpha|}} \left( \begin{array}{c}
|\alpha|\\
\alpha \\
\end{array} \right)
p^\alpha s_\alpha.$$

Again, it is well known that ${\mathcal F}(\mu)=0$ implies $\mu=0$.

\subsection{Cauchy transform} Among the various choices of Cauchy's integral transforms
we prefer the following one.  The {\it Cauchy transform} of a measure $\mu \in \cM(\R^d)$ is
$$\mathcal C \mu (z) = \int_{\R^d} \frac{d\mu(x)}{(x_1-z_1)...(x_d-z_n)},\ \ z \in \C^n \setminus \R^d.$$
This is an analytic function in the variable $z$, which determines
the measure $\mu$.

\begin{proposition} If $\mu \in \cM(\R^d)$ and
$\mathcal C\mu = 0$, then $\mu = 0$.
\end{proposition}

\begin{proof} Let
$\hat{\mu}(u)  = \int_{\R^d} e^{-i u\cdot x} d\mu(x)$ be the
Fourier transform of $\mu$. Due to the decay assumptions,
$\hat{\mu} \in {\mathcal C}^\infty (\R^d)$ and $\hat{\mu}$ is
uniformly bounded on $\R^d$.

Let $y \in \R^d$ be a vector with positive entries. Then Laplace's
transform is well defined, and Fubini's Theorem yields:
$$ \int_{\R_+^d} e^{-u\cdot y} \hat{\mu}(u) du =  \int_{\R_+^d} e^{-u\cdot y} \int_{\R^d} e^{-i u\cdot x} d\mu(x)\ du=
$$ $$ \int_{\R^d} \int_{\R_+^d} e^{-u \cdot (y+ix)} du\  d\mu(x) = \int_{\R^d} \frac{d\mu(x)}{(y_1 + ix_1)...(y_d + i x_d)}
= $$ $$ (-i)^n \mathcal C \mu( -iy)  = 0.$$

Therefore $\hat{\mu} (u) = 0$ for all $u \in \R_+^d$. By repeating
the same computations to $u$ in a different octant, with the
proper choice of the signs of $y_k$, we find $\hat{\mu} = 0$,
whence $\mu=0$.

\end{proof}

The characterization of all Cauchy transforms is available via
standard techniques of harmonic analysis in a product of half
planes, cf. \cite{Koranyi,SW}. We briefly indicate below the
details.

\begin{thm} An analytic function $F(z_1,...,z_d), \ z_i \notin
\R$ for all $i=1,...,n$, satisfies $F = {\mathcal C} \mu$ for a
positive, finite measure $\mu$ on $\R^d$ if and only if the
functions $f_k(z)$ defined inductively as
$$ f_1(z_1,...,z_d) =
(F(z_1,...,z_{d-1},z_d)-F(z_1,...,z_{d-1},\overline{z_d}))/(2i),$$
$$ f_{k+1}(z) = $$ $$(f_k(z_1,...,z_{k-1}, z_k, z_{k+1},...,z_d)-
f_k(z_1,...,z_{k-1}, \overline{z_k}, z_{k+1},...,z_d))/(2i),$$ for
all $1 \leq k \leq d,$ fulfill the conditions:
\begin{enumerate}
\item $ f_d(z) \geq 0,\ \ \ \Im z_i >0, \ 1 \leq i \leq d,$
\item $ \sup_y \int_{\R^d} |f_d(x_1+iy_1,...,x_d+iy_d)|
dx_1...dx_d
< \infty,$
\item $\lim_{y_k \rightarrow \infty} F(z_1,...,z_{k-1},x_k+iy_k,
z_{k+1},...,z_d) = 0, \  1 \leq k \leq d.$
\end{enumerate}

In this case $w^\ast-\lim_{y \downarrow 0} \frac{f_n(x+iy)
dx_1...dx_d}{\pi^d} = \mu$.
\end{thm}

The case $n=1$ corresponds to the well studied class of Cauchy
transforms of positive measures on the line, see for instance
\cite{A}.

\begin{proof} We outline the main steps of the proof in the case
of two variables $(z,w)$. The passage to higher dimensions does
not require additional arguments.

Suppose that $F={\mathcal C}\mu$. The recursively generated
function is in this situation an iterated Poisson transform:
$$ f_2(x+iy,u+iv) = \int_{\R^2}
\frac{y}{(x-s)^2+y^2}\ \frac{v}{(u-t)^2+v^2} d\mu(s,t).$$ Thus
$f_2$ is non-negative for $y,v>0$ and
$$ \int f_2(x+iy,u+iv) dx du = \pi^2 \int_{\R^2} d\mu(s,t)
<\infty.$$ Moreover, $w^\ast-lim_{v \downarrow 0} \frac{v du
}{u^2+v^2} = \pi \delta(u)$, whence the stated
$w^\ast$-convergence of $f_2$ holds true.

Conversely, assume that the analytic function $F(z,w)$ satisfies
the conditions in the statement. The function $f_2$ is then
harmonic in every variable and by a repeated application of
Theorem II.2.5 of \cite{SW} there exists a finite, positive
measure $\mu$ on $\R^2$ such that $f_2$ is its iterated Poisson
transform. Thus, for all $z,w \notin \R$ we have:
$$ F(z,w)-F(z,\overline{w})-F(\overline{z},w)+F(\overline{z},
\overline{w}) = $$ $$ \int_{\R^2}
[\frac{1}{t-z}-\frac{1}{t-\overline{z}}][\frac{1}{s-w}-\frac{1}{s-\overline{w}}]
d\mu(s,t).$$ By polarization, we obtain the identity
$$F(z,w)-F(z,w')-F(z',w)+F(z',w') =
$$ $$ \int_{\R^2}
[\frac{1}{t-z}-\frac{1}{t-z'}][\frac{1}{s-w}-\frac{1}{s-w'}]
d\mu(s,t).$$

By letting $z' \rightarrow \infty$ and then $w' \rightarrow
\infty$ we obtain
$$ F(z,w) = \int_{\R^2}
\frac{d\mu(s,t)}{(t-z)(s-w)}.$$
\end{proof}

\subsection{Determinateness and quasi-analyticity} Without aiming
at obtaining optimal results (they will be improved later via
Hilbert space methods) we exemplify below the link between
quasi-analyticity of a Fourier-Laplace transform of a measure
supported by the positive octant in the Euclidean space and its
determinateness. In this direction we go back to some old theorems
due to Bochner-Taylor and Bochner \cite{BT,Bochner}.

We state a result of \cite{BT} for the convenience of the reader.
This is only a first generic result. Quite a few ramifications of
it are presented in the same article as well as in
\cite{Bochner,Ehrn,Hormander}.

\begin{thm} Let $f$ be a $\mathcal C^\infty$
function defined on a domain $\Omega$ of $\R^d$. Let $x_0 \in
\Omega$ and let $m_n$ be a sequence of positive numbers. Assume
that
\begin{enumerate}
\item  $ \sum_{|\alpha|=n} |\frac{\partial^\alpha}{\partial x^\alpha} f (x)|^2 \leq m_n^2$
for all $n \geq 0$ and $x \in \Omega$;
\item  $\frac{\partial^\alpha}{\partial x^\alpha} f(x_0) = 0$ for all $\alpha \in \N^d$;
\item  $\sum_n m_n^{-1/n} = \infty$.
\end{enumerate}
Then $f$ is identically zero on $\Omega$.
\end{thm}

The application to moment problems is now transparent:

\begin{thm} Let $\{s_\alpha\}$ be the moment sequence of a measure
$ \mu \in \cM(\R_+^d)$. If
$$ \sum_{k=0}^\infty [\sum_{|\alpha|=k} |s_\alpha|^2]^{-\frac{1}{2k}} = \infty,$$
then the moment problem is determinate on $\R_+^d$.
\end{thm}

\begin{proof} Denote by $\mu_1, \mu_2$ two measures having the prescribed
moment sequence
$\{s_\alpha;\alpha \in \N_0^d\}$. Let $\mu = \mu_1-\mu_2$ and take
the Laplace transform $\tilde{\mu}.$ This is a $\mathcal C^\infty$
function on the closed octant $\R_+^d$, and by assumption
$\frac{\partial^\alpha}{\partial x^\alpha} \tilde{\mu} (0) = 0$ for all $\alpha \in \N_0^d$.
Note that
$$ |\frac{\partial^\alpha}{\partial x^\alpha} \tilde{\mu}(x)| \leq 2 s_\alpha, \ \  \alpha \in \N_0^d.$$

If $n=1$, then $\tilde{\mu}$ extends smoothly to the negative
semi-axis and the statement follows from a standard
quasi-analyticity criteria \cite{Bochner}. We prove the general
statement by induction on $d$.

Due to the induction hypothesis, Laplace transform $\tilde{\mu}$
vanishes with all partial derivatives along the boundary of
$\R_+^n$. Whence we can extend $\mathcal C^\infty$
 it by zero to the whole $\R^n$. Then one applies directly Theorem 1 of
 \cite{BT}.
 \end{proof}
  The article \cite{BT} and its companion \cite{Bochner} contain similar criteria of quasi-analyticity,
  in terms of the growth of $L(x,\partial)^k  \tilde{\mu}$, where $L(x,\partial)$ is a linear partial
  differential operator. The multivariate quasi-analyticity theme is also amply discussed in
  \cite{Ehrn}. A link between the rate of rational approximation and quasi-analyticity is
  developed in \cite{Plesniak}. It is worth mentioning that all these authors obtain
  their quasi-analyticity conditions via Bernstein's theorem on the real axis, by
  restriction to lines or
  one-dimensional curves.

  Next we give an example of a variation of the above result,
  derived from Theorem 10 of \cite{BT} with the same proof as above.

  \begin{corollary}  If $\{s_\alpha\}$ is the moment sequence of a measure $\mu \in \cM(\R_+^d)$ such
  that $$ \sum_{k=0}^\infty [\sum_{|\alpha|=k} \left( \begin{array}{c} k\\
             \alpha \end{array} \right) s_{2\alpha}]^{-\frac{1}{k}} = \infty,$$
             then $\mu$ is determinate on $\R_+^d$.
             \end{corollary}

\section{Hilbert space methods}
The spectral theory of self-adjoint operators in Hilbert space is
well suited and provides powerful techniques for the study of the
the moment problem and in particular of the determinacy question.
The present section is devoted to this approach. A thorough treatment of the  one dimensional moment problem in terms of  the self-adjoint extension theory of symmetric operators can be found in the recent survey by B. Simon \cite{Simon}.

Let $L$ be a linear a linear functional on the $*$-algebra
$\C[x]\equiv \C[x_1,\dots,x_d]$. Throughout this section we assume
that $L$ is {\it positive} functional, that is, $L(p\bar{p}) \geq
0$ for all $p \in \C[x]$ or equivalently $L(q^2)\geq 0$ for all
$q\in \R[x]$.

The technical tool that relates the moment problem to operator
theory in Hilbert space is the {\it GNS-construction}. We briefly
develop this construction. (Note that it works for positive linear
functionals on arbitrary unital $*$-algebras.) Since $L$ is
positive, the Cauchy-Schwarz inequality
$$
|L(p_1\overline{p_2})|^2 \leq
L(p_1\overline{p_1})L(p_2\overline{p_2}),~~ p_1,p_2 \in \C[x],
$$
holds and implies that $\cN_L:=\{ p\in \C[x]: L(p\overline{p})=0
\}$ is an ideal of the algebra $\C[x]$. Hence there exist a scalar
product $\langle \cdot,\cdot\rangle_L$ on the quotient space
$\cD_L=\C[x]/ \cN_L$ and an algebra homomorphism $\pi_L$ of
$\C[x]$ into the linear operators on $\cD_L$ defined by
$$\langle
p+\cN_L,q+\cN_L\rangle_L =L(p\overline{q}) ~{\rm and}~
\pi_L(p)(q+\cN_L)=pq+\cN_L,~~p,q \in \C[x].
$$
 Let
$\Hh_L$ denote the Hilbert space completion of the pre-Hilbert
space $\cD_L$. If no confusion can arise we omit all subscripts
$L$ and write $q$ for $q+\cN_L$ and $X_k$ for $\pi(x_k)$,
$k{=}1,\dots,d$. Then we have $\pi(p)q=pq$ and
\begin{align}\label{Lrep}
\langle \pi(p)p_1,p_2\rangle = \langle p_1,\pi(\overline{p})p_2
\rangle =L(pp_1\overline{p_2}),~~ p,p_1,p_2 \in \C[x].
\end{align}
In particular,  $L(p)=\langle \pi(p)1,1 \rangle $ for $p \in
\C[x]$ and  $X_1,\dots,X_d$ are {\it pairwise commuting symmetric
operators on the dense invariant domain} $\cD=\pi(\C[x])1$ of the
Hilbert space $\Hh$.

The next proposition relates the moment problem to spectral
measures of strongly commuting $d$-tuples of self-adjoint
operators extending the $d$-tuple $(X_1,\dots,X_d)$. This
connection goes back to the early days of functional analysis.
\begin{proposition}\label{mp-spec}
A positive linear functional $L$ of $\C[x_1,\dots,x_d]$ is a moment functional if and only if there exists a $d$-tuple $(A_1,\dots,A_d)$ of strongly commuting self-adjoint operators $A_1,\dots,A_d$ acting on a Hilbert space $\cK$ such that $\Hh$ is a subspace of $\cK$ and $X_1\subseteq A_1,\dots, X_d \subseteq A_d$.\\
If $L$ is a moment functional, then all solutions $\mu$ of the
moment problem for $L$ are of the form $\mu(\cdot)= \langle
E_{(A_1,\dots,A_d)}(\cdot)1,1 \rangle $, where $(A_1,\dots,A_d)$
is such a $d$-tuple and $E_{(A_1,\dots,A_d)}$ denotes its spectral
measure.
\end{proposition}

Let us first recall the notions occurring in this proposition. A
$d$-tuple $(A_1,\dots,A_d)$ of self-adjoint operators on a Hilbert
space $\cK$ is called {\it strongly commuting} if the resolvents
$(A_k{-}iI)^{-1}$ and $(A_l{-}iI)^{-1}$ commute or equivalently if
the corresponding spectral measures $E_{A_k}$ and $E_{A_l}$
commute (that is,  $E_{A_k}(M)E_{A_l}(N)=E_{A_l}(N)E_{A_k}(M)$ for
all Borel subsets $M,N$ of $\R^d$) for all $k,l=1,\dots,d, k{\neq}
l$. Such a $d$-tuple has a unique spectral measure
$E_{(A_1,\dots,A_d)}$ on $\R^d$ determined by
\begin{align}\label{jointspec}
E_{(A_1,\dots,A_d)}(M_1\times\cdots \times M_d)
=E_{A_1}(M_1)\cdots E_{A_d}(M_d)
\end{align}
for arbitrary Borel subsets $M_1,\dots,M_d$ of $\R$.\\

\noindent Proof of Proposition \ref{mp-spec}:

Let $L$ be a moment functional for a measure $\mu \in \cM(\R^d)$.
Then it is easily checked that the multiplication operators $A_k$,
$k=1,\dots,d$, by the coordinate functions $x_k$ form a $d$-tuple
of strongly commuting self-adjoint operators on the Hilbert space
$\cK=L^2(\mu)$ such that $\Hh \subseteq \cK$ and $X_k\subseteq
A_k$ for $k=1,\dots,d$.
 The spectral measure $E\equiv E_{(A_1,\dots,A_d)}$ of this $d$-tuple is given by $E(M)f=\chi_M \cdot f$, $f \in L^2(\mu)$, where $\chi_M$ is the characteristic function of the Borel set $M\subseteq \R^d$. Hence we  have $\langle E(M)1,1 \rangle =\mu(M)$.

Conversely, let $(A_1,\dots,A_d)$ be such a $d$-tuple and let
$E_{(A_1,\dots,A_d)}$ be its spectral measure. Put
$\mu(\cdot)=\langle E(\cdot)1,1 \rangle$. Let $p \in \C[x]$. Since
$X_k \subseteq A_k$, we have $p(X_1,\dots,X_d) \subseteq
p(A_1,\dots,A_d)$. The polynomial $1$ belongs to the domain of
$p(X_1,\dots,X_d)$ and hence of $p(A_1,\dots,A_d)$. From the
spectral calculus we obtain
\begin{align*}
\int p(\lambda)~ d\mu(\lambda) &= \int p(\lambda)~ d\langle
E_{(A_1,\dots,A_d)}(\lambda)1,1 \rangle  =\langle p(A_1,\dots,A_d) 1,1 \rangle  \\&
= \langle p(X_1,\dots,X_d)1,1 \rangle  = L(p(x_1,\dots,x_d)),
\end{align*}
where the last equality follows from (\ref{Lrep}). This shows that
$\mu$ is a solution of the moment problem for $L$. \hfill $\Box$

\medskip
\noindent For the determinacy problem we have the following
important sufficient criterion which was first noticed by A. Devinatz \cite{Devinatz}, see also \cite{Fug}.
\begin{proposition}\label{op-det}
Suppose  $L$ is a moment functional and $\mu$ is a representing
measure for $L$. If each symmetric operator $X_k$, $k=1,\dots,d$,
is essentially self-adjoint (that is, $\overline{X_k}=X_k^*$) or
equivalently if $\C[x_1,\dots,x_d]$ is dense in
$L^2((1+x_k^2)\mu)$ for each $k=1,\dots,d$, then the moment
problem for $L$ is determinate and the polynomials
$\C[x_1,\dots,x_d]$ are dense in $L^2(\mu)$.
\end{proposition}

\noindent Our proof of Proposition \ref{op-det} given below is
different from the ones  in \cite{Devinatz} and \cite{Fug}. It is essentially based on
the following  lemma.
\begin{lemma}\label{self-sym}
Let $A$ be a closed symmetric operator on a Hilbert space $\cK$
and let $P$ denote the orthogonal projection of $\cK$ onto its
closed subspace $\Hh$. Suppose that $\cD$ is a dense linear
subspace of $\Hh$ such that $\cD\subseteq \cD(A)$, $A\cD \in \Hh$
and $X:=A\lceil \cD$ is an essentially self-adjoint operator on
$\Hh$. Then we have $PA\subseteq AP$. Moreover, if $Y$ denotes the
restriction of $A$ to $(I-P)\cD(A)$, then $A= \overline{X}\oplus
Y$ on $\cK=\Hh \oplus \Hh^\perp$.
\end{lemma}
Proof. Suppose that $\varphi \in \cD(A)$. Let $\psi \in \cD$.
Using the facts hat  $\psi$ and $X\psi$ are in $\Hh$ and the relation
$X\subseteq A$, it follows that
$$
\langle X\psi,P\varphi \rangle =\langle PX\psi,\varphi  \rangle =
\langle X\psi,\varphi \rangle=$$ $$ \langle A\psi,\varphi \rangle
= \langle \psi,A\varphi \rangle = \langle P\psi,A\varphi \rangle
=\langle \psi,PA\varphi \rangle.
$$
Since $\psi \in \cD$ was arbitrary, $P\varphi \in \cD(X^\ast)$ and
$X^\ast P\varphi =AP\varphi$. Since $X$ is essentially
self-adjoint, $X^\ast=\overline{X}$ and so $AP\varphi
=\overline{X} P\varphi = PA\varphi$. That is, $PA\subseteq AP$.

If $\varphi \in \cD(A)$, we have $P\varphi \in \cD(\overline{X})$
and $(I-P)\varphi \in \cD(A)$ as just shown. Hence $\varphi
=P\varphi +(I-P)\varphi \in \cD(\overline{X}) \oplus \cD(Y)$ and
$A\varphi = \overline{X}\varphi + Y(I-P)\varphi$, so $A \subseteq
\overline{X}\oplus Y$. By definition the converse inclusion is
clear. \hfill $\Box$

\medskip
\noindent Proof of Proposition \ref{op-det}:

 Suppose that the operators $X_1,\dots,X_d$ are essentially self-adjoint on $\Hh$.

 Let $(A_1,\dots,A_d)$ be a $d$-tuple of strongly commuting self-adjoint operators on a larger Hilbert space $\cK$ such that $X_j\subseteq A_j$, $j{=}1,\dots,d$. By Lemma \ref{self-sym} there is a decomposition $A_j =X_j \oplus Y_j$ on $\cK=\Hh \oplus \Hh^\perp$. Obviously, $Y_j$ is also self-adjoint and  $E_{A_j}=E_{\overline{X_j}}\oplus E_{Y_j}$ for the corresponding spectral measures. Hence the spectral measures of $\overline{X_k}$ and $\overline{X_l}$ commute and  formula (\ref{jointspec}) yields
$E_{(\overline{X_1},\dots,\overline{X_d})} \subseteq
E_{(A_1,\dots,A_d)}$. Therefore, $\mu(\cdot):= \langle
E_{(\overline{X_1},\dots,\overline{X_d})}(\cdot)1,1 \rangle$ is
equal to $ \langle E_{(A_1,\dots,A_d)}(\cdot)1,1 \rangle$. Since
each  solution of the moment problem is of the form $\langle
E_{(A_1,\dots,A_d)}(\cdot)1,1 \rangle$ by Proposition
\ref{mp-spec} and so coincides with $\mu$ by the preceding, the
moment problem is determinate.

Let $M$ be a Borel set of $\R^d$. By definition the polynomials are dense in $\Hh_L$, so there is a sequence $\{p_n\}$ of polynomials $p_n \in \C[x]$ such that $p_n(x) \to
E_{(\overline{X_1},\dots,\overline{X_d})}(M)1$ in $\Hh_L$ and or equivalently
$$
||(p_n(x)- E_{(\overline{X_1},\dots,\overline{X_d})}(M))1||^2 = \int_{\R^d} |p_n(\lambda)- \chi_M(\lambda)|^2 d\langle E_{(\overline{X_1},\dots,\overline{X_d})}(\lambda)1,1\rangle \to 0
$$
as $n \to \infty$. Hence the closure of polynomials $\C[x]$ in $L^2(\mu)$ contains all characteristic functions $\chi_M$, so it is equal to $\Hh_L=L^2(\mu)$.

Since $X_k$ is essentially self-adjoint, $(X_k\pm i)\C[x]$ is dense in $\Hh_L=L^2(\mu)$ which obviously implies the density of $\C[x]$ in $L^2((1+x_k^2)\mu)$.

Conversely, if $\C[x]$ is dense in $L^2((1+x_k^2)\mu)$, it follows that $(x_k\pm i)\C[x]=(X_k\pm i)\C[x]$ is dense in $L^2(\mu)$ and so in its subspace $\Hh_L$. Hence $X_k$ is essentially self-adjoint on the Hilbert space $\Hh_L$.
\hfill $\Box$

\medskip
The preceding Proposition \ref{op-det} suggests the
following definitions.
\begin{definition}
Let $\mu \in \cM(\R^d)$ and let $L$ be its moment functional, that
is, $L(p) =\int p~d\mu$ for $p \in \C[x]$. We say that $\mu$, respectively
$L_\mu$, is {\rm strongly determinate} if the symmetric operators
$X_1,\dots,X_d$ are essentially self-adjoint on the Hilbert space
$\Hh_{L_\mu}$. We call $\mu$ resp. $L_\mu$ {\rm strictly
determinate} if $\mu$ resp. $L_\mu$ is determinate and the
polynomials $\C[x_1,\dots,x_d]$ are dense in $L^2(\mu)$.
\end{definition}

Using these notions we can restate Proposition \ref{op-det} as
follows: The measure $\mu$ is strongly determinate if and only if
$\C[x_1,\dots,x_d]$ is dense in $L^2((1+x_k^2)\mu)$ for
$k=1,\dots,d$. Strong determinacy always implies strict
determinacy and so determinacy.

Strong determinacy and another concept called ultradeterminacy
have been introduced  by B. Fuglede, see \cite{Fug}, p. 57. A
measure $\mu \in \cM(\R^d)$ resp. its moment functional $L_\mu$ is
called {\it ultradeterminate} if the polynomials
$\C[x_1,\dots,x_d]$ are dense in $L^2((1+||x||^2)\mu)$. Since the
norm of $L^2((1+||x||^2)\mu)$ is obviously stronger than the norm
of $L^2((1+x_k^2)\mu)$ for  $k=1,\dots,d$, each ultradeterminate
measure is strongly determinate.

From the Weierstrass approximation theorem it follows that each
measure with compact support is ultradeterminate. Further, if $\mu
\in \cM(\R^d)$ and $\C[x_1,\dots,x_d]$ is dense in $L^p(\mu)$ for
some $p> 2$, then $\mu$ is ultradeterminate, see \cite{Fug}, p.61.

The equivalence (i)$\leftrightarrow$(iii) of Propositon
\ref{det-1} below shows that in the one-dimensional case all four
concepts ultradeterminacy, strong determinacy, strict determiancy
and determinacy are the same. However, in dimension $d\geq 2$,
they are all different. Examples of a strongly determinate measure
that is not ultradeterminate and of a strictly determinate measure
that is not stronlgy determinate have been given in \cite{S1}. By
Theorem 5.4 in \cite{Bergth} there exist (rotation invariant)
determinate measures on $\R^d$, $d\geq 2$, which are not strictly
determinate.

Examples show that ultradeterminacy and of strong determinacy are
rather strong concepts. In many situations they are even too
strong. It seems to us that strict determinacy is a more important
and fundamental concept. One reason is that strict determinacy is
crucial for the reduction theorem \ref{applth} proved below.
Another reason stems from the theory of orthogonal polynomials,
because for a strict determinate measure $\mu$ any sequence of
orthonormal polynomials is an orthonormal basis of $L^2(\mu)$.

Let us briefly discuss the preceding Propositions \ref{mp-spec}
and \ref{op-det} in the one-dimensional case $d=1$. From
Proposition \ref{mp-spec} it follows at once that a positive
functional on $\C[x_1]$ is always a moment functional, because
each symmetric operator has a self-adjoint extension in a larger
Hilbert space (see e.g. \cite{AG}, Nr. 111). A more careful
analysis (see \cite{A} or \cite{Simon}) shows that the operator
$X_1$ is either essentially self-adjoint or has deficiency indices
$(1,1)$. For the determinacy in dimension one we have the
following fundamental result.

\begin{proposition}\label{det-1}
For a measure $\mu\in \cM(\R)$ and its moment functional $L_\mu$
the following statements are equivalent:
\begin{itemize}
\item{\em (i)} $\mu$ is determinate.
\item{\em (ii)} The symmetric operator $X_1$ is essentially self-adjoint on $\Hh_{L_\mu}$.
\item{\em (iii)} $\C[x_1]$ is dense in $L^2((1+x_1^2)\mu)$.
\item{\em (iv)} $1$ is in the closure of $(x_1{+}z)\C[x_1]$ in $\Hh_{L_\mu}$ for some (resp. all) $z \in \C\setminus \R$.
\item{\em (v)} $1$ is in the closure of $(x_1{+}z)\C[x_1]$ in $L^2(\mu)$ for some (resp. all) $z \in \C\setminus \R$.
\end{itemize}
\end{proposition}
\noindent Proof. See \cite{A} Chapter 3, or \cite{Simon}. \hfill
$\Box$

\medskip
The striking difference between the one-dimensional moment problem
and the multi-dimensional case $d \geq 2$ appears also in the
Hilbert space approach. First, there exist positive linear
functionals on $\C[x_1,\dots,x_d]$, $d\geq 2$, which are not
moment functionals. Probably the simplest example of this kind is
Example 6 in \cite{S2}. To explain the difference to the
one-dimensional case, let $L$ be a positive linear functional on
$\C[x_1,\dots,,x_d]$, $d \geq 2$. Since the complex conjugation on
$\C[x_1,\dots,x_d]$ commutes with the operator $X_k$, it is a
bijective map of ${\ker}~ (X_k^*{-}zI)$ on ${\ker}~
(X_k^*{-}\bar{z}I)$. Therefore, each symmetric operator $X_k$ has
equal deficiency indices and so a self-adjoint extension on the
Hilbert space $\Hh_L$. However, in order to apply Proposition
\ref{mp-spec} one needs {\it strongly commuting} self-adjoint
extensions of $X_1,\dots,X_d$ in a possible larger Hilbert space.
To find reasonable results for the existence of such  extensions
is a very difficult task. This makes the Hilbert space approach
much more complicated and so less powerful than in the
one-dimensional case. Secondly, as noted above there are various
different determinacy concepts in case $d \geq 2$. In particular,
we don't know of any useful general {\it necessary}
operator-theoretic criterion for determinacy in the
multi-dimensional case.

\section{Polynomial and rational approximation}

In this section we develop a number of criteria for determinacy
which are useful by their flexibility in applications. All of them
are based on density conditions of polynomials.

Let us begin with two simple results concerning the following
general question: Given $\mu \in \cM(\R^d)$, what can be said
about the set $V_\mu$ of equivalent measures ?

The next proposition is a classical result of M.A. Naimark
\cite{Nai} for the one-dimensional case. The proof given in
\cite{A}, p. 47, carries over verbatim to the multi-dimensional
case.
\begin{proposition}\label{naimark}
Suppose that $\mu \in \cM(\R^d)$ and $\nu \in V_\mu$. Then
$\C[x_1,\dots,x_d]$ is dense in $L^1(\nu)$ if and only if $\nu$ is
an extreme point of the convex set $V_\mu$. In particular,
$\C[x_1,\dots,x_d]$ is dense in $L^1(\mu)$ if $\mu$ is
determinate.
\end{proposition}

\begin{proposition}
Suppose that $\nu, \mu \in \cM(\R^d)$ and $\nu \cong \mu$.
\begin{itemize}
\item{\em (i)} If $C\subseteq \R^d$ is a real
 algebraic variety and ${\rm supp}~\mu \subseteq C$, then ${\rm supp}~\nu \subseteq C$.
\item{\em (ii)} If $p \in \C[x_1,\dots,x_d]$ is bounded on ${\rm supp}~\mu$, then $p$ is also
bounded on ${\rm supp}~ \nu$ and we have ${\sup}~ \{ |p(x)|;x \in
{\rm supp}~ \mu \} = {\sup}~ \{|p(x)|;x \in {\rm supp}~ \nu \}$.
\end{itemize}
\end{proposition}
Proof. (i): By definition there are polynomials $p_1,\dots,p_r \in \R[x]$ such
that $C{=}\{x\in \R^d:p_1(x)=\cdots p_r(x)=0\}$. Since $\nu \cong
\mu$ and $C \subseteq {\rm supp}~\mu$, we have
$$
-\int p_j(x)^2 ~d\nu(x)= -\int p_j^2(x) ~d\mu(x) =-\int_{C}
p_j(x)^2~d\mu(x)=0
$$
which implies that $p_j(x)=0$ on ${\rm supp}~\nu$,
$j=1,\dots,r$. Therefore, ${\rm supp}~\nu \subseteq C$.

(ii): In this proof we repeat an argument which has been used in
\cite{S2}, p. 228. Let $K_\nu$ and $K_\mu$ denote the supremum of
$|p(x)|$ over ther sets ${\rm supp}~ \nu$ and ${\rm supp}~ \mu$,
respectively. For $\lambda >K_\mu$, put $M_\lambda :=\{x \in
\R^d:|p(x)| \geq \lambda \}$. Then  we obtain
$$
\lambda^{2n} \nu(M_\lambda) \leq \int_{M_\lambda} (p\bar{p})^n
~d\nu \leq \int_{\R^d} (p\bar{p})^n ~d\nu = \int_{\R^d}
(p\bar{p})^n ~d\mu \leq K_\mu^{2n} \mu(\R^d)
$$
for  $n \in \bbN$. Since $\lambda > K_\mu$, the preceding
inequality can be only valid for arbitrary $n$ when
$\nu(M_\lambda)=0$. Therefore, ${\rm supp}~\nu \subseteq \{x\in
\R^d:|p(x)|\leq K_\mu\}$ which in turn implies that $K_\nu \leq
K_\mu$. Interchanging the role of $\nu$ and $\mu$, we get
$K_\nu=K_\mu$. \hfill $\Box$

\medskip
Now we return to the determinacy problem.

For a Borel mapping $\varphi:\R^d\to \R^m$ and a Borel measure
$\mu$ on $\R^d$, we denote by $\varphi(\mu)$ the image of $\mu$
under the mapping $\varphi$, that is, $\varphi(\mu)(M):=
\mu(\varphi^{-1}(M))$ for any Borel set $M$ of $\R^m$. Then the
transformation formula
\begin{align}\label{immeasure}
\int_{\R^m} f(y)~ d\varphi(\mu)(y) = \int_{\R^d} f(\varphi(x))~
d\mu(x)
\end{align}
holds for any function $f \in L^1(\varphi(\mu))$.

Let $\pi_j(x_1,\dots,x_d)=x_j$ be the j-th coordinate mapping of
$\R^d$ into $\R$. Then $\pi_j(\mu)$ is the j-th marginal measure
of $\mu$. The following basic result was obtained by L. C.
Petersen \cite{Petersen}.
\begin{thm}\label{pet}
If all marginal measures $\pi_1(\mu),\dots,\pi_d(\mu)$ are
determinate, then $\mu$ is determinate.
\end{thm}

In order to make our exposition complete in the multidimensional case as promised above, we repeat the proof of Theorem \ref{pet} from \cite{Petersen}.

\smallskip
Proof. Suppose that $\nu \in V_\mu$. Let $\chi_1,\dots,\chi_d$ be characteristic functions of Borel subsets of $\R$. If $p_1,\dots,p_d$ are polynomials in one variable, we estimate
\begin{align*}
||\chi_1&(x_1)\dots \chi_d(x_d)-p_1(x_1)\dots p_d(x_d)||_{L^1(\nu)} \\
\leq &||(\chi_1(x_1)-p_1(x_1))\chi_2(x_2)\dots \chi_d(x_d)||_{L^1(\nu)}+\\
& ||p_1(x_1)(\chi_2(x_2)-p_2(x_2))\chi_3(x_3)\dots \chi_d(x_d)||_{L^1(\nu)}
+\dots+ \\
&||p_1(x_1)\dots p_{d-1}(x_{d-1})(\chi_d(x_d)-p_d(x_d))||_{L^1(\nu)}\\
\leq & ||\chi_1(x_1)-p_1(x_1)||_{L^2(\nu)}~ ||\chi_2(x_2)\dots \chi_d(x_d)||_{L^2(\nu)} +\\ &||\chi_2(x_2)-p_2(x_2)||_{L^2(\nu)}~||p_1(x_1)\chi_3(x_3)\dots \chi_d(x_d)||_{L^2(\nu)}+\dots+ \\ &
||p_1(x_1)\dots p_{d-1}(x_{d-1})||_{L^2(\nu)}~ ||\chi_d(x_d)-p_d(x_d)||_{L^2(\nu)}.
\end{align*}
Let $j\in \{1,\dots,d\}$. Clearly, $\nu \cong \mu$ implies that
$\pi_j(\nu)\cong \pi_j(\mu)$, so that $\pi_j(\nu)=\pi_j(\mu)$, because $\pi_j(\mu)$
is determinate by assumption. Hence we have
\begin{align}\label{l2mu}
||\chi_j(x_j){-}p_j(x_j)||_{L^2(\nu)}=||\chi_j(x_j){-}p_j(x_j)||_{L^2(\pi_j(\nu))}
=\end{align} \begin{align*}
||\chi_j(x_j){-}p_j(x_j)||_{L^2(\pi_j(\mu))} .
\end{align*}
Since $\pi_j(\mu)$ is determinate, the polynomials $\C[x_j]$ are
dense in $L^2(\pi_j(\mu))$. Therefore, by (\ref{l2mu}), we can
first choose $p_1$ such that
$||\chi_1(x_1){-}p_1(x_1)||_{L^2(\nu)}$ becomes arbitrarily
small, then $p_2$ such that
$||\chi_2(x_2){-}p_2(x_2)||_{L^2(\nu)}$ becomes small and finally
$p_d$ such that $||\chi_d(x_d){-}p_d(x_d)||_{L^2(\nu)}$ becomes
small. By the above inequality, $||\chi_1(x_1)\dots
\chi_d(x_d)-p_1(x_1)\dots p_d(x_d)||_{L^1(\nu)}$ becomes as small
as we want. Since the span of functions $\chi_1(x_1)\dots
\chi_d(x_d)$ is dense in $L^1(\nu)$, this shows that the
polynomials are dense in $L^1(\nu)$, so $\nu$ is an extreme point
of   $V_\mu$ by Proposition \ref{naimark}. Because $\nu \in V_\mu$
was arbitrary, $\nu =\mu$.
 \hfill $\Box$

\smallskip
As shown by a simple example in \cite{Petersen},
there exist
determinate measures for which not all marginal measures are determinate.

Now we derive two theorems which are corollaries of Proposition
\ref{pet}. They provide some  general determinacy criteria in
terms of polynomial approximations and contain, of course,
Theorem \ref{pet} as special cases.

Let $\varphi=(\varphi_1,\dots,\varphi_m):\R^d\to \R^m$ be a
polynomial mapping, where $\varphi_1,\dots,\varphi_m$ $\in
\R[x_1,\dots,x_d]$. Since $\mu \in \cM(\R^d)$, we have
$\varphi(\mu) \in \cM(\R^m)$ and $\varphi_k(\mu) \in \cM(\R)$ by
(\ref{immeasure}). Combining Proposition \ref{det-1} and formula
(\ref{immeasure}) it follows that the following three conditions
are equivalent:
\begin{itemize}
\item{\em (a)} $\varphi_k(\mu)$ is determinate.
\item{\em (b)} There exist a sequence $\{q_{kn};n{\in}\bbN\}$ of polynomials
$q_{kn}\in \C[y_1]$ such that
$(\varphi_k(x){+}z)q_{kn}(\varphi_k(x)) \to 1$ as $n\to \infty$ in
$\Hh_{L_\mu}$ for some $z \in \C\setminus \R$.
\item{\em (c)} $(\varphi_k(x){+}z)\C[\varphi_k(x)]$ is dense in $L^2(\varphi_k(\mu))$ for some $z \in \C\setminus \R$.
\end{itemize}
Conditions (b) and (c) can be rephrased by saying that $1$ resp.
elements of $L^2(\varphi_k(\mu))$ can be approximated by (certain)
polynomials. It should be emphasized that condition (b) is
"intrinsic" to the moment problem, that is, it depends only on the
moment functional $L_\mu$, but not on the particular representing
measure $\mu$.
\begin{thm}\label{pedcor1}
Let $K$ be a closed subset of $\R^d$ and $\mu \in \cM(K)$. Suppose
that the map  $\varphi:K\to \R^m$ is injective and that for each
$k=1,\dots,m$ one of the equivalent conditions (a)--(c) is
satisfied. Then $\mu$ is determinate on $K$.
\end{thm}
Proof. Let $\nu$ be another measure such that $\nu \cong \mu$ and
${\rm supp}~\nu \subseteq K$. Since $\varphi$ is a polynomial
mapping and $\nu \cong \mu$, it follows from (\ref{immeasure})
that $\varphi(\nu) \cong \varphi (\mu)$. By assumption, all
marginal measures $\pi_k(\varphi(\mu))=\varphi_k(\mu)$ are
determinate. Therefore, $\varphi(\mu)$ is determinate by
Proposition \ref{pet}, so that $\varphi(\nu)=\varphi(\mu)$. Since
${\rm supp}~\varphi(\mu)\subseteq \varphi(K)$ and  $\varphi:K\to
\R^m$ is injective by assumption, we conclude that $\nu =\mu$.
\hfill $\Box$.

\smallskip
\noindent For the second result we consider $r$ polynomial
mappings $\varphi^j :\R^d \to \R^{m_j}$, $j{=}1,\dots,r$, and
define a polynomial mapping
$\varphi=(\varphi^1,\dots,\varphi^r):\R^d \to \R^m$, where
$m=m_1+\cdots+m_r$.
\begin{thm}\label{pedcor2}
Let $\mu$ and $K$ be as in Theorem \ref{pedcor1}. If
$\varphi:K\to \R^m$ is injective and all measures
$\varphi^1(\mu),\dots, \varphi^r(\mu)$ are strictly determinate,
then $\mu$ is determinate on $K$.
\end{thm}
Proof. The proof is based on a generalization of Proposition
\ref{pet} which is obtained by replacing the map $\pi_j$ onto a
single coordinate by a map $\pi^j$ onto a finite set of
coordinates. To be more precise, we write $y \in \R^m$ as
$y=(y_{11},\dots,y_{1m_1},y_{21},\dots,y_{rmr})$ and define
mappings $\pi^j:\R^d\to \R^{m_j}$ by
$\pi^j(y)=(y_{j1},\dots,y_{jm_j})$.  Then we have the following
result: {\it If $\nu \in \cM(\R^m)$ and all measures
$\pi^1(\nu),\dots,\pi^r(\nu)$ are strictly determinate , then
$\nu$ is determinate.}

A proof of this statement can be given by repeating the arguments
of the above proof of Proposition \ref{pet} and using that the
polynomials $\C[y_{j1},\dots,y_{jm_j}]$ are dense in
$L^2(\pi^j(\nu))$ for the strictly determinate measure
$\pi^j(\nu)$. We do not carry out the details.

Using this generalization instead of Proposition \ref{pet} we can
argue as in the proof of the preceding theorem. \hfill $\Box$

\smallskip
\noindent The presence of "sufficiently many" bounded polynomials on the set
${\rm supp}~\mu$ allows us to prove stronger results that the
plain determinacy. As a sample we consider subsets of cylinders
with compact base sets.
\begin{proposition}
Let $K$ be a closed subset of  $\R^d, d\geq 2$, such that $K$ is a subset
of $K_0 \times \R$, where $K_0$ is a compact set of
$\R^{d-1}$. Suppose that $\mu \in \cM(K)$. If the marginal measure
$\pi_d(\mu)$ is determinate, then $\mu$ is ultradeterminate and
hence strongly determinate.
\end{proposition}
Proof. We shall write $x \in \R^d$ as $x=(y,x_d)$ with $y\in
\R^{d-1}$ and $x_d \in \R$. Suppose that $f(y)$ and $g(x_d)$ are
given continuous functions with compact support. Let us abbreviate
$M_1={\rm sup} ~\{|g(x_d)|;x_d \in \R \}$, $M_2={\rm sup}~\{
|p(y)|;y \in K_0\}$ and $M_3={\sup}~\{||y||;y \in K \}$. Moreover
we  denote by  $||\cdot||_1$ the norm of $L^2((1+||x||^2)\mu)$, by
$||\cdot||_2$ the norm of $L^2(\mu)$ and by $||\cdot||_3$ the norm
of $L^2((1+x_d^2)\pi_d(\mu))$. Let $p \in \C[y]$ and $q \in
\C[x_d]$. Using the assumption ${\rm supp}~ \mu \subseteq K
\subseteq  K_0 \times \R$ and formula (\ref{immeasure}), we
estimate
\begin{align*}
&||f(y)g(x_d)-p(y)q(x_d)||_1 \leq
||(f(y)-p(y))g(x_d)||_1 + ||p(y)(g(x_d)-q(x_d))||_1\\
&\leq ||f(y)-p(y)||~M_1 + M_2M_3 ||(g(x_d)-q(x_d))(1+x_d^2)||_2\\
& \leq~{\rm sup}~ |f(y)-p(y)|;y \in K_0\} ~M_1 \mu(\R^d)+ M_2M_3
||g(x_d)-q(x_d)||_3.
\end{align*}
By the Weierstrass approximation theorem, $p \in \C[y]$ can be
chosen such that the supremum of $|f(y){-}p(y)|$ over the compact
set $K_0$ is arbitrary small. Since the marginal measure
$\pi_d(\mu)$ is determinate, $\C[x_d]$ is dense in
$L^2((1+x_d^2)\pi_d(\mu))$ by Proposition \ref{det-1}. Hence we
can choose $q \in \C[x_d]$ such that $||g-q||_3$ is sufficiently
small. Therefore, we have shown that the function $f(y)g(x_d)$ is
in the closure of $\C[x]$ in $L^2((1+||x||^2)\mu)$. Since the span
of such functions is obviously dense, it follows that the
polynomials $\C[x]$ are dense in  $L^2((1+||x||^2)\mu)$, so $\mu$
is ultradeterminate. \hfill $\Box$\\

We reproduce below from \cite{PS} another general determinateness
criterion related to polynomial approximation.

\begin{proposition}\label{suffdetermined}
Let $K$ be a closed subset of $\R^d$, and let $\mu, \nu \in
\mathcal{M}(K).$ Let $f\in \mathcal C(K,\R)$ be a function
satisfying $f\ge1$, a.e.\ on $K$ (with respect to both $\mu$ and
$\nu$).

Assume that there exists a sequence of polynomials $p_n$ in
$\R[x]$ such that $p_n\to\frac1f$ under the norm
$||\cdot||=||\cdot||_{L^2(\mu)} +||\cdot||_{L^2(\nu)}$, and let
$$\mathcal A_0:= \mathcal A_0(K,f):=\Bigl\{\frac p{f^k}\colon p\in\R[x],\ k\ge0,\
\frac{p(x)}{f(x)^k}\to0\text{ for }|x|\to\infty,\ x\in K\Bigr\}.$$
If $\mathcal A_0$ separates the points of $K$, then $\mu=\nu$.
\end{proposition}

\begin{proof}
All fractions $\frac p{f^k}$ (with $p\in\R[x]$ and $k\ge0$) are
integrable with respect to both $\mu$ and $\nu$. We prove by
induction on $k$ that
$$\int\frac p{f^k}d\mu=\int\frac p{f^k}d\nu.$$
For $k=0$ this holds, since $L_\mu=L_\nu$ by assumption. For the induction step $k\to
k+1$ note that
$$\frac p{f^{k+1}}=\lim_{n\to\infty}\Bigl(p_n\cdot\frac p{f^k}
\Bigr),$$ both in $L^1(\mu)$ and  $L^1(\nu)$ (since
$||fg||_1^2=\langle|f|,|g| \rangle^2\le||f||_2^2\cdot||g||_2^2$ by
the Cauchy-Schwarz inequality), and hence
$$\int\frac{p}{f^{k+1}}\>d\mu=\int\frac{p}{f^{k+1}}\>d\nu.$$
Let $\mathcal B$ denote the subalgebra of $\mathcal C(K,\R)$
consisting of the functions $\phi$ for which the limit of
$\phi(x)$, for $x\in K$ and $|x|\to\infty$, exists in $\R$. So
$\mathcal B= \mathcal C(K^+,\R)$ where $K^+=K\cup\{\infty\}$ is
the one-point compactification of $K$. (If $K$ is already compact,
put $K^+=K$.) By assumption, the subalgebra $\mathcal A:=\R1\oplus
\mathcal A_0$ of $\mathcal B$ separates the points of $K^+$.
Therefore, by the Stone-Weierstra\ss\ theorem, $\mathcal A$ is
dense in $\mathcal B$ under uniform convergence. Therefore the
measures $\mu$ and $\nu$ coincide as linear functionals on
$\mathcal B$, because if $b\in \mathcal B$ and $a_n$ in $\mathcal
A$ with $a_n\to b$ under $||\cdot||_ \infty$, then
$\int_Ka_n\to\int_Kb$ for both $\mu$ and $\nu$ since
$\mu(K)=\nu(K)<\infty$. It follows that $\mu=\nu$.
\end{proof}

The proposition applies to a variety of situations. For
instance:\\
\begin{itemize}

 \item {\bf 1.}\ Assume we
have $f\in\R[x]$ with $f\ge1$ on $K$ for which there is a sequence
$\{ p_n \}$ in $\R[x]$ with $||1-fp_n||_{L,2}\to0$. If $\mathcal
A_0(K,f)$ separates the points of $K$, then the moment problem on
$K$ given by $L$ is determinate. (Indeed, $||\frac1f-p_n||_2\le
||1-fp_n||_2$ since $f\ge1$.) Note that, contrary to the preceding
proposition, the condition $||1-fp_n||_{L,2}\to0$ is intrinsic in
$L$ and its values on polynomials.\\

\item {\bf 2.} \ Assume that an entire function is given in form of its
Taylor series
$$\frac1f=\sum_\alpha c_\alpha x^\alpha.$$
Assume that the algebra $\mathcal A_0(K,f)$ fulfills the
separation condition in the statement of the Proposition, and let
$\mu$ be a positive measure on $K$ with moments $\{ s_\alpha\}$.
The normal convergence condition
$$\limsup_\alpha\>\bigl(|c_\alpha|\cdot s_{2\alpha}\bigr)^
{1/|\alpha|}<1$$ will assure that the partial sums converge to
$1/f$ in $L^2(\mu)$, and hence the above result is applicable.\\

\item {\bf 3.}\ There are natural choices of continuous functions $f\ge1$
for which $\mathcal A_0(K,f)$ separates the points of $K$, like
$f=1+\sum_ix_i^2$ or $f=e^{x^2}$.\\

\item {\bf 4.}\ In the one-dimensional case we know  from Proposition \ref{det-1} that the converse of the assertion of Proposition \ref{suffdetermined} holds as well: The moment functional $L$ is  determinate if and only if there exists a sequence of
polynomials $\{ p_n\}$ for which $||1-(1+x^2)p_n||_{\Hh_L}\to 0$.
\end{itemize}

\section{Partial determinacy, moment functionals and determinacy}

In this section we show how partial determinacy can be used to conclude that a positive linear functional is even a moment functional and to ensure determinacy. The pioneering  result on partial determinacy is due A. Devinatz \cite{Devinatz}. A slightly stronger theorem was obtained later by G.I. Eskin \cite{Eskin}, see e.g. \cite{Nussbaum} and the references in \cite{Eskin}.

 For technical reasons we first develop all results in the case $d=2$. The case $d \geq 3$ is more complicated and will be treated afterwards. The following theorem is closely related to Theorem 22 in \cite{Fried}.
\begin{thm}\label{detmf}
Let $\cQ \subseteq \C[x_1,x_2]$ be a set of polynomials such that $\C[x_1,x_2]$ is the linear span of polynomials $p(x_1)q(x_1,x_2)$, where $p \in \C[x_1]$ and $q\in \cQ$. Suppose that $L$ is a positive linear functional on $\C[x_1,x_2]$ such that for each  $q \in \cQ$ the (positive) linear functional $L_q(p):=  L(p(x_1)(x_2^2+1)q\overline{q})$ on $\C[x_1]$ is determinate. Then $L$ is a moment functional.
\end{thm}

The special case $\cQ{=}\{x_2^k;k \in \N_0\}$ of Theorem
\ref{detmf} is  due to Eskin \cite{Eskin}. It is the first
assertion of the following corollary.
\begin{corollary}\label{Eskin}
Suppose  $s=\{s_{n,k};k,n \in \N_0\}$ is a positive semidefinite $2$-sequence such that
for each $k{\in} \N_0$ the $1$-sequence
\begin{align}\label{eskin}
\{s_{(n,2(k+1))}+s_{(n,2k)}; n \in \N_0 \}
\end{align}
 is determinate. Then  $s$ is a moment sequence. If in addition the sequence $\{s_{0,n};n \in \N_0\}$ is determinate, then $s$ is determinate.
\end{corollary}
 Proof. Put $\cQ{=}\{x_2^k; k \in \N_0\}$. Since $s_{(n,m)}=L(x_1^nx_2^m)$, $L_{x_2^k}$ has the moment sequence (\ref{eskin}), so $L$ is a moment functional by Theorem \ref{detmf}.

 The determinacy of the sequence $\{s_{(n,2)}+s_{(n,0)}\}$ clearly implies that of the sequence $\{s_{(n,0)}\}$. Hence both marginal sequences of $s$ are determinate, so $s$ is determinate by Petersen's Theorem \ref{pet}. \hfill $\Box$

\medskip
\noindent
The main technical parts of the proof of Theorem \ref{detmf} are contained in the following two lemmas. If $A$ is a self-adjoint operator, the unitary operator  $U_A:=(A-iI)(A+iI)^{-1}$ is called the {\it Cayley transform} of $A$.
\begin{lemma}\label{extab}
Let  $A_j, j \in J,$ be a family of pairwise strongly commuting self-adjoint operators $A_j$ and let $B$ be a densely defined symmetric operator acting  on a Hilbert space $\Hh$. Suppose that $U_{A_j}BU_{A_j}^*=B$ for all $j \in J$. Then there exists a family $\tilde{A_j}, j \in J,$ of strongly commuting self-adjoint operators $\tilde{{A_j}}$ and a self-adjoint operator $\tilde{B}$ on a Hilbert space $\cK \supseteq \Hh$ such that $\tilde{B}\supseteq B$, $A_j$ and $B$ strongly commute and $\tilde{A_j}\supseteq {A_j}$ for all $j \in J$.
\end{lemma}

Proof. The proof is based on the so-called "doubling trick" which is often used in operator theory.  Define $\tilde{{A_j}}={A_j} \oplus {A_j}$ for $j \in J$ and $B_0=B \oplus(-B)$ on the Hilbert space $\cK=\Hh\oplus \Hh$. Obviously, the operators $\tilde{A_j},j \in J,$ strongly commute as the operators $A_j,j \in J,$ do by assumption. In the proof we freely use the self-adjoint extension theory of symmetric operators (see \cite{AG} or \cite{RN}). Without loss of generality we can asume that the operator $B$ is closed. Let $\cN_\pm (T)={\rm ker}(T^*\pm iI)$ denote the deficiency spaces of a symmetric operator $T$. From the definition of $B_0$ we obtain
$\cN_\pm (B_0)= \cN_\pm (B)\oplus \cN_\mp (B).$
Hence $V(\varphi,\psi):=(\psi,\varphi)$, where $\varphi \in \cN_+(B), \psi \in \cN_-(B)$,  defines an isometric linear map $V$ of $\cN_+(B_0)$ onto $\cN_-(B_0)$. The restriction $\tilde{B}$ of $(B_0)^*$ to the domain $\cD(\tilde{B}):=\cD(B_0) + (I{-}V)\cN_+(B_0)$ is a self-adjoint extension of $B_0$.

From the assumption $U_{A_j}BU_{A_j}^*=B$ we derive  $U_{A_j}B^*U_{A_j}^*=B^*$ which in turn yields $U_{A_j}\cN_\pm(B)=\cN_\pm(B).$ From  $\cN_\pm (B_0)= \cN_\pm (B)\oplus \cN_\mp (B)$ and the definition of $V$ it follows that the Cayley transform $U_{\tilde{{A_j}}}=U_{A_j}\oplus U_{A_j}$ of  $\tilde{{A_j}}$
maps $(I{-}V)\cN_+(B_0)$
onto itself. Since $U_{A_j}^*\cD(B)=\cD(B)$ by assumption and hence $U_{A_j}\cD(B)=\cD(B)$, we get $U_{\tilde{{A_j}}}\cD(\tilde{B})=\cD(\tilde{B}).$ Combined with the relation  $$U_{\tilde{{A_j}}}(B_0)^*U_{\tilde{{A_j}}}^* =U_{A_j}B^*U_{A_j}^*\oplus (-U_{A_j}B^*U_{A_j}^*)=B^*\oplus (-B^*) = (B_0)^*$$ the latter implies that $U_{\tilde{{A_j}}}\tilde{B}U_{\tilde{{A_j}}}^*=\tilde{B}$ and so $U_{\tilde{{A_j}}}U_{\tilde{B}}U_{\tilde{{A_j}}}^*=U_{\tilde{B}}$. Hence the resolvents of $\tilde{{A_j}}$ and $\tilde{B}$ commute. Therefore,   $\tilde{{A_j}}$ and $\tilde{B}$ strongly commute.
\hfill $\Box$

\begin{lemma}\label{extab1}
Let $A$ and $B$ be closed symmetric linear operators ona Hilbert space $\Hh$ and let $\cD\subseteq \cD(A)\cap\cD(B)$ be a  linear subspace  such that $A\cD \subseteq \cD$, $B\cD \subseteq \cD$, and $AB\psi=BA\psi$ for $\psi \in \cD$. Let $\cQ$ be a subset of $\cD$ such that the linear span of vectors $p(A)\varphi$, where $p \in \C[x]$ and $\varphi \in \cQ$, is dense in $\Hh$ and a core for $B$ . For $\varphi \in \cQ$ we denote by $A_\varphi$ the restriction of $A$ to the invariant linear subspace  $\cD_\varphi:=\C[A](B+i)\varphi$.  Suppose that for each $\varphi \in \cQ$ the symmetric operator $A_\varphi$ is essentially self-adjoint in the Hilbert space obtained as the closure of $\cD_\varphi$ in $\Hh$. Then the operator $A$ on $\Hh$ is self-adjoint and we have $U_ABU_A^* =B$.
\end{lemma}

Proof. Let $\varphi \in \cQ$, $z \in \C\backslash \R$ and $p \in \C[x]$. Since $A_\varphi$ is essentially self-adjoint, there exists a sequence $\{p_n\}$ of polynomials $p_n \in \C[x]$ such that $(A+z)p_n(A)(B+i)\varphi\to p(A)(B+i)\varphi$ in $\Hh_\varphi$. Using that $B$ is a symmetric operator commuting with $A$ on $\cD$ we compute
\begin{align*}
&||(A+z)p_n(A)(B+i)\varphi-p(A)(B+i)\varphi||^2 =\\
& ||(B+i)((A+z)p_n(A)\varphi-p(A)\varphi)||^2 =\\
&||B((A+z)p_n(A)\varphi-p(A)\varphi) ||^2 + ||(A+z)p_n(A)\varphi-p(A)\varphi ||^2 =\\ &||(A+z)Bp_n(A)\varphi-Bp(A)\varphi ||^2 + ||(A+z)p_n(A)\varphi-p(A)\varphi ||^2
\end{align*}
From these equations we conclude that $(A+z)Bp_n(A)\varphi \to Bp(A)\varphi$ and $(A+z)p_n(A)\varphi  \to p(A)\varphi$ as $n \to \infty$. Since the  span of vectors  $p(A)\varphi$ is dense in $\Hh$, the preceding shows that the range of $(A+z)$ is dense in $\Hh$, so the  operator  $A$ on $\Hh$ is self-adjoint. Further, because $(A+z)^{-1}$ is bounded, it follows that $p_n(A)\varphi\to (A+z)^{-1}p(A)\varphi$ and $Bp_n(A)\varphi \to (A+z)^{-1}Bp(A)\varphi$ as $n\to \infty$. Since the operator $B$ is closed, the latter yields  $B(A+z)^{-1}p(A)\varphi =(A+z)^{-1}Bp(A)\varphi$. By assumption the span of vectors $p(A)\varphi$ is a core for $B$. Therefore, we obtain that $B(A+z)^{-1}\psi =(A+z)^{-1}B\psi$ for all vectors $\psi \in \cD(B)$. Setting $z=i$ and $z=-i$, it follows that $U_A=I-2i(A+i)^{-1}$ and $U_A^{-1}=I+2i(A-i)^{-1}$ map the domain $\cD(B)$ into itself. Therefore, we conclude that $U_A\cD(B)=\cD(B)$. Hence the relation
$$
U_AB\psi=(I-2i(A+i)^{-1})B\psi=B(I-2i(A-i)^{-1})\psi=BU_A\psi
$$
for $\psi \in \cD(B)$ implies that $U_ABU_A^* =B$.
\hfill $\Box$

\medskip
 Proof of Theorem \ref{detmf}. We shall apply Lemma \ref{extab1}. Let $A$ and $B$ be the closures of the multiplication operators  by $x_1$ and $x_2$, respectively, on the Hilbert space $(\Hh_L,\langle\cdot, \cdot \rangle)$. For $q \in \cQ$, we denote by $(\Hh_q,\langle\cdot, \cdot \rangle_q)$ the Hilbert space of the moment functional $L_q$ on $\C[x_1]$. From the definitions of the functional $L_q$ and of the scalar products (see formula (\ref{Lrep})) we obtain
$$
\langle p_1(x_1),p_2(x_1)\rangle_q=
L(p_1(x_1)\overline{p_2}(x_1)(x_2^2+1)q\overline{q})
= \langle p_1(A)(B+i)q,p_2(A)(B+i)q \rangle
$$
for $p_1,p_2 \in \C[x_1]$. Since $L_q$ is determinate by assumption, the multiplication operator $X_1$ is essentially self-adjoint on $\Hh_q$ by Proposition \ref{det-1},(i)$\to$(ii). Therefore, by the preceding equality the operator $A_{q}$ is essentially self-adjoint on $\cD_{q}$ for each $q \in \cQ$. (Here we retain the notation of Lemma \ref{extab1}.) Hence the assumptions of Lemma \ref{extab1} are fulfilled, so we have  $U_ABU_A^* =B$. Therefore, combining the conclusion of Lemma \ref{extab}, applied to the single operator $A$,  with Proposition \ref{mp-spec} we infer that $L$ is a moment functional.
 \hfill $\Box$

\medskip
\noindent
Now we turn to the case when the dimension $d$ is larger than $2$. The following theorem is the counter-part of
Corollary \ref{Eskin}.

\begin{thm}\label{eskinmulti}
Let $s=\{s_n;n \in \N_0^d\}$ be a positive semidefinite
$d$-sequence, where $d \geq 3$. Suppose that for arbitrary numbers
$j=1,\dots,d{-}1$ and $k_1,\dots$, $k_{j-1},k_{j+1},\dots,k_{d}
\in 2 \N_0$ the $1$-sequence
\begin{align}\label{det1}
\{s_{(k_1,\dots,k_{j-1},n,k_{j+1},\dots,k_{d-1},k_d+2)}
+s_{(k_1,\dots,k_{j-1},n,k_{j+1},\dots,k_{d-1},k_d)}; n \in
\N_0\},
\end{align}
is determinate. Further, suppose that for all numbers $j,l\in\
\{1,\dots,d{-}1\}$, $j<l$, all sequences of one of the following
two sets of  $1$-sequences indexed over $n \in \N_0$
\begin{align}\label{det2}
 \{s_{(k_1,\dots,k_{j-1},n,k_{j+1},\dots,k_{l-1},k_l+2,k_{l+1}, \dots,k_d)} +
 s_{(k_1,\dots,k_{j-1},n,k_{j+1},\dots,k_{l-1},k_l,k_{l+1},\dots,k_d)} \},
\end{align}
\begin{align}\label{det3}
\{s_{(k_1,\dots,k_{j-1},k_j+2,k_{j+1},\dots,k_{l-1},n,k_{l+1},\dots,k_d)}
+
s_{(k_1,\dots,k_{j-1},k_j,k_{j+1},\dots,k_{l-1},n,k_{l+1},\dots,k_d)}\},
\end{align}
where $k_1,\dots,k_d \in 2 \N_0$, are determinate. Then the $d$-sequence $s$ is a moment sequence.
If in addition the $1$-sequence
$\{s_{(0,\dots,0,n)};n \in \N_0\}$ is determinate, then $s$ is also determinate.
\end{thm}

Proof. The proof is given by some modifications of the above proof of Theorem \ref{detmf}. Let $A_j$, $j=1,\dots,d$, denote the closure of the multiplication operator by $x_j$ in the Hilbert space $\Hh_{L_s}$. First we fix $j\in \{1,\dots,d{-}1\}$ and apply Lemma \ref{extab1} in the case where $A=A_j$, $B=A_d$, and  $\cQ$ is the set of all monomials $x_1^{k_1}\cdots x_{j-1}^{k_{j-1}}x_{j+1}^{k_{j+1}}\cdots x_d^{k_d}$. Since all sequences
 (\ref{det1}) are determinate, the assumptions of Lemma \ref{extab1} are fulfilled. Therefore, by Lemma \ref{extab1}, the operator $A_j$ is self-adjoint and we have $U_{A_j}A_dU_{A_j}^* = A_d$. In order to show that the self-adjoint operators $A_1,\dots,A_{d-1}$ strongly commute we apply Lemma \ref{extab1} once more. Let $j,l\in \{1,\dots,d{-}1\}$, $j<l$. Assume that all sequences (\ref{det2})  are determinate. Then we set $A=A_j$, $B=A_l$,  and let $\cQ$ be the set of monomials $x_1^{k_1}\cdots x_{j-1}^{k_{j-1}}x_{j+1}^{k_{j+1}}\cdots x_d^{k_d}$. From Lemma \ref{extab1} we obtain that $U_{A_j}A_lU_{A_j}^* = A_l$. Since $A_j$ and $A_l$ are both self-adjoint, the latter yields  $U_{A_j}U_{A_l}U_{A_j}^* = U_{A_l}$. Hence the resolvents of $A_j$ and $A_l$ commute which implies that $A_j$ and $A_l$ strongly commute. If the sequences (\ref{det3}) of the second set are determinate, we interchange the role of $j$ and $l$ and proceed in a similar manner. Thus, $A_1,\dots,A_{d-1}$ is a family of strongly commuting self-adjoint operators such that $U_{A_j}A_dU_{A_j}^* = A_d$ for $j=1,\dots,d{-}1$. That is, the assumptions of Lemma  \ref{extab} are satisfied with $B:=A_d$. From Lemma \ref{extab} and Proposition \ref{mp-spec} we conclude that $s$ is a moment sequence. As in proof of Corollary \ref{eskin} the determinacy assertion follows from Petersen's result (Theorem \ref{pet}).
\hfill $\Box$

\section{Quasi-analyticity}

Growth conditions on the moments are another kind of assumptions
that yield determinacy criteria for the moment problem. However,
it is well known, even in the one variable case, that growth
conditions cannot characterize determinateness. In this short
section we will elaborate this method in the case of
quasi-analyticity.

If $T$ is a symmetric operator on a Hilbert space,  a vector $\varphi \in \cap_{n=1}^\infty  \cD(T^n)$  is called a {\it quasi-analytic vector for $T$} if
\begin{align}\label{quasivec}
\sum_{n=1}^\infty  ||T^n\varphi ||^{-1/n} =+\infty.
\end{align}
In terms of the moments $s_n:=\langle T^n \varphi,\varphi \rangle
$,  (\ref{quasivec}) is obviously equivalent to the so-called {\it Carleman
condition}
\begin{align}
\label{carle}
\sum_{n=1}^\infty  s_{2n}^{-1/{2n}} =+\infty.
\end{align}

The importance of this condition for the moment problem stems form
the following classical result due to T. Carleman \cite{Carleman}
in dimension one. See also \cite{Hormander}.

\begin{thm}\label{Car} If $\{s_n;n\in \N_0\}$ is a moment
sequence satisfying Carleman' s condition (\ref{carle}), then the corresponding moment problem is determinate.
             \end{thm}

Carleman's result is the main technical ingredient of
the following theorem.

\begin{thm}\label{qavector}
Suppose that $s=\{s_n;n \in \N_0^d\}$ is a positive semi-definite $d$-sequence such that the first $d{-}1$ marginal $1$-sequences
\begin{align}\label{det4}
\{s_{(n,0,\dots,0))};n \in \N_0\}, \{s_{(0,n,\dots,0)};n \in \N_0\} , \dots, \{s_{(0,\dots,0,n,0)};n \in \N_0\}
\end{align}
 satisfies the Carleman condition (\ref{carle}). Then $s$ is a moment sequence.\\ If in addition the last marginal $1$-sequence $\{s_{(0,\dots,0,n)};n \in \N_0\}$ satisfies  (\ref{carle}), then the moment sequence $s$ is also determinate.
\end{thm}
 Proof.  That the sequences in (\ref{det4}) satisfy condition (\ref{carle}) means that the vector $1$ of the corresponding Hilbert space $\Hh_L$ is quasi-analytic for the symmetric operators $X_1,\dots, X_{d-1}$. Assume first that $d=2$. Since the operator $(X_2+i)X_2^k$ commutes with $X_1$, the vector $(X_2+i)X_2^k 1$ is also quasi-analytic for $X_1$ according to Theorem 4 in \cite{Nussbaum}. Using (\ref{Lrep}) we see that the moment $1$-sequence for $X_1$ corresponding to this vector is just the sequence in  (\ref{eskin}), so it fulfills (\ref{carle}) as well and  is determinate by Carleman's Theorem \ref{Car}. Therefore, $s$ is a moment sequence by Corollary \ref{Eskin}.

In the case $d\geq 3$ the proof is similar.  Using again the fact that commuting operators preserve quasi-analytic vectors
it follows that all sequences in equations  (\ref{det1}), (\ref{det2}) and (\ref{det3}) satisfy Carleman's condition, so they are determinate by Theorem \ref{Car}. Hence $s$ is a moment sequence by Theorem \ref{eskinmulti}.

If in addition (\ref{carle}) holds for the sequence $\{s_{(0,\dots,0,n)}\}$, then $s$ is determinate by Petersen's Theorem \ref{pet}.
\hfill $\Box$

\smallskip
\noindent
The preceding theorem contains the following fundamental result due to A.E. Nussbaum \cite{Nussbaum}: {\it A positive semi-definite $d$-sequence is
a determinate moment sequence if all $d$ marginal $1$-sequences satisfy Carleman's condition.}

Other classical papers on quasi-analytic or analytic vectors are \cite{Chernoff}, \cite{Chernoff1}, \cite{KM}, \cite{Ne} and \cite{Nussbaum1}.

\section{Towards a parametrization of solutions}

The parametrization of all solutions of an indeterminate one
variable moment problem on the line or on the semi-axis is well
understood from at least two perspectives: via a canonical
linear-fractionar representation of the Cauchy transforms of the
respective measures, or via an analysis of the generalized
self-adjoint extensions of a symmetric operators. See for details
\cite{A,Buchwalter,Simon,ST}.

The aim of the present section is to parametrize all solutions to
an indeterminate moment problem in the positive octant of $\R^d$,
by "integrating" a family of Stieltjes 1D moment problems with
respect to a pencil of lines spanning the octant. Although rather
theoretical, at this point, our approach combines for the first
time Nevanlinna's parametrization with the characterization of
Fantappi\`e transforms of positive measures.

We recall first some known facts about Nevalinna parametrization of
solutions to the Hamburger and Stieltjes moment problems, following the conventions in
\cite{Simon}. Let $\mu$ be a rapidly decaying at infinity
positive measure on the real line. For every $z \notin \R$ and
every equivalent measure $\nu \in V_\mu$, that is in our notation
$\nu \cong \mu$, the values of the Cauchy transform
$$ \int_{\R} \frac{d\nu(x)}{x-z}$$
fill a disk $D(\mu,z)$, known as the Weyl disk. Denoting by $P_k(x), k \geq
0,$ the sequence of orthogonal polynomials associated to the
measure $\mu$ (or any equivalent measure, as these polynomials
depend solely on the sequence of moments) and by $Q_k$ the
associated orthogonal polynomials of the second kind:
$$ Q_k(x) = \int_{\R} \frac{P_k(x)-P_k(u)}{x-u} d\mu(u), \ k \geq
0,$$ the equation of the boundary of $D(\mu,z)$, in the case $\Im
z >0$ is
$$ D(\mu,z) = \{ w \in \C; \ \frac{w-\overline{w}}{z-\overline{z}}
= \sum_{n=0}^\infty |w P_k(z)+Q_k(z)|^2\}.$$ The radius of this
circle is
$$ \rho(z) = \frac{1}{|z- \overline{z}|}
\frac{1}{|\sum_{n=0}^\infty |P_n(z)|^2}$$
and $\rho(z_0) >0$ for an arbitrary non-real $z_0$, if and only if the
moment problem associated to $\mu$ is indeterminate. If this is the case,
then
$$ A(z) = z \sum_{k=0}^\infty Q_k(0)Q_k(z),$$
$$ B(z) = -1 + z \sum_{k=0}^\infty Q_k(0)P_k(z),$$
$$ C(z) = 1+ z\sum_{k=0}^\infty P_k(0)Q_k(z),$$
$$ D(z) = z \sum_{k=0}^\infty P_k(0)P_k(z),$$
are entire functions fully constructible from the moments of $\mu$.

The convex set $V_\mu$ is then parametrized by all analytic maps
$\Phi$ from the upper half-plane to the closure of the upper
half-plane in the Riemann sphere, as follows:
$$  \int_{\R} \frac{d\nu(x)}{x-z} = - \frac{C(z) \Phi(z)+ A(z)}{D(z)\Phi(z)+B(z)}.$$
The functions $\Phi$, known as Nevanlinna, or Herglotz functions, provide
a free parameter into the bijective description of all $\nu \in V_\mu$.

The case of an indeterminate Stieltjes problem requires additional care.
Specifically, assume that $\mu$ is a positive measure on $\R_+$, having all moments finite.
The orthogonal polynomials $P_k, Q_k, k \geq 0$, as well as the entire functions
$A,B,C,D$ are the same as above. It is only the parametrization of all measures
$\nu \in V_\mu(\R_+)$ which is more restrictive, in the following precise sense.
There exists a bijective correspondence between the set $\nu \in V_\mu(\R_+)$
and all Nevanlinna functions of the form
$$\Psi(z) = d_0 + \int_0^\infty \frac{d\sigma(t)}{t-z}, \ \ \Im z >0,$$
where
$$ d_0 = - \lim_{k \rightarrow \infty} \frac{Q_k(0)}{P_k(0)},$$
and $\sigma$ is an arbitrary positive measure, rapidly decaying at infinity.
The formula relating the measure $\nu$ with the function $\Psi$ is the same
as in the case of Hamburger problem:
$$  \int_{\R_+} \frac{d\nu(x)}{x-z} = - \frac{C(z) \Psi(z)+ A(z)}{D(z)\Psi(z)+B(z)}.$$
Equivalently, the functions $\Psi$ entering into the above representation
can be characterized by the fact that either they are constant, or they are
analytic in $\C \setminus \R_+$, and in any case
$$ \Psi(t) \geq d_0, \ \ {\rm for \ all}\ t<0.$$
For details see \cite{Simon}.

All these known facts being recalled, we return now to the
multivariate moment problem in the positive octant of $\R^d$,
where we will combine the above structure and the theory of the
Fantappi\`e transform. To this aim, start with a measure $\mu \in
\cM(\R_+^d)$. For every vector $a \in {\rm int}\R_+^d$ we consider
the push-forward measure $(a\cdot )\mu = (a \cdot )_\ast \mu$
defined by the identity
$$ \int_{\R_+}  \phi(t) d[(a \cdot )_\ast \mu](t) = \int_{{R_+^n}}
\phi (a \cdot x) d\mu(x).$$ For an equivalent measure $\nu \in
V_\mu(\R_+^d)$ the moments of these push-forward measures agree,
that is $(a \cdot ) \nu \in V_\mu(\R_+)$. Therefore, for a fixed
$z \notin \R_+$ we have
$$ \int_{{\R_+^d}} \frac{d\mu(x)}{a\cdot x - z} =
\int_{\R_+} \frac{d[(a \cdot ) \nu](t)}{t-z} \in \Delta((a \cdot )
\nu, z)$$ and these Cauchy transforms depend on a
"free" parameter $\Psi$ as explained above.

Next we reverse this process, and analyze the characteristic
properties of a family of positive measures $\sigma_a \in
\cM(\R_+)$ indexed over $a \in {\rm int}\R_+^d$ which satisfy
$\sigma_a \in V_{(a \cdot)\mu}(\R_+)$ and  "integrate" to a
measure $\nu \in V_\mu(\R_+^d)$. First of all these measures
should satisfy a smoothness dependence on the parameter $a$.
Second, they should fulfill a compatibility condition with respect
to the homothety $a \mapsto \lambda a$ where $\lambda>0$:
\begin{equation}
\sigma_{\lambda a} = \lambda_\ast \sigma_a, \ \ a \in {\rm
int}\R_+^d, \ \lambda >0.
\end{equation}
In other terms, for every bounded continuous functions $\phi$ we
have
$$ \int_{{\R_+}} \phi(t)d\sigma_{\lambda a}(t) =
 \int_{{\R_+}} \phi(\lambda t)d\sigma_{a}(t).$$
To such a family of measures $\sigma_a$ we associate the function
$$ F(a,z) = \int_{\R_+} \frac{d\sigma_a(t)}{t-z},\ \ a \in {\rm
int}\R_+^d, \ z \notin \R_+.$$ Notice that for every $\lambda>0$
$$ F(\lambda a , \lambda z) = \int_{\R_+} \frac{d\sigma_{\lambda a}(t)}{t-\lambda
z}= \int_{\R_+} \frac{d\sigma_a(t)}{\lambda t- \lambda z} =
\lambda^{-1} F(a,z).$$ Hence, in view of the characterization of
Fanttapi\`e transforms of positive measures (see Theorem
\ref{Fantappie}) we are led to the following rather abstract
observation.

\begin{thm} Let $\mu \in \mathcal M(\R_+^d)$. For every $a \in {\rm
int}\R_+^d$ let $P_k(a,z), Q_k(a,z)$ be the orthogonal polynomials
of the first and second kind associated to the push forward
measure $(a\cdot )_\ast \mu$ on $\R_+$. Assume that all these
measures $(a\cdot )_\ast \mu, \  a \in {\rm int}\R_+^d$ are
indeterminate, so that the associated entire functions $A(a,z),
B(a,z),$ $C(a,z), D(a,z)$ and limit $d_0(a)$ exist.

Then there exists a bijective correspondence between all measures
$\nu \in V_\mu(\R_+^d)$ and the $(-1)$-homogeneous, $\mathcal
C^\infty$ cross sections in the Weyl disks bundle of the form
$$ F(a,z) \in \Delta((a\cdot )_\ast \mu,z),$$
$$ F(a,z)  = - \frac{C(a,z) \Psi_a(z)+ A(a,z)}{D(a,z)\Psi_a(z)+B(a,z)}, \ z \in \C \setminus \R_+, $$
where $\Psi_a$ are Nevanlinna functions analytic in $\C \setminus
\R_+$ satisfying $\Psi_a(t) \geq d_0(a), \  t<0,$ and such that
the function $F(a,-t), \ \ a \in {\rm int}\R_+^d, t>0$ is
completely monotonic.

\end{thm}

\begin{proof} One direction of the correspondence was established
before. In the other direction, the assumptions in the statement imply, in virtue of Nevanlinna parametrization theorem, the existence of a positive measure $\sigma_a$ on $\R_+$
with the property
$$ F(a,z) = \int_{\R_+} \frac{d\sigma_a(t)}{t-z}.$$
Moreover, the dependence of $\sigma_a$ on the parameter $a$ is
smooth, as stated. In virtue of Theorem \ref{Fantappie}, there
exists a positive measure $\nu \in \cM(\R_+^d)$ satisfying
$$ F(a,-t) = \int_{{\R_+^d}} \frac{d\nu(x)}{a\cdot x +t}.$$
On the other hand, for all $a$ and $t$,
$$ F(a,-t)  =
\int_{\R_+}\frac{d\sigma_a(u)}{u+t}.$$ Consequently the measures
$(a\cdot )\nu$ and $\sigma_a$ have the same one dimensional
moments, that is $(a\cdot )\nu$ and $(a\cdot )\mu$ have the same
moments. In particular, for every fixed vector $a$ and positive
integer $k$ we obtain
$$ \int (a\cdot x)^k d\mu = \int (a\cdot x)^k d\nu.$$
A Vandermonde determinant argument implies then
$$ \int x^\alpha d\mu = \int x^\alpha d\nu, \ \ |\alpha|=k.$$
In other terms $\nu \in V_\mu(\R_+^d).$
\end{proof}

A similar analysis can be carried out for an indeterminate moment
problem in $\R^d$ provided that a characterization of Fantappi\`e
transforms of positive measures $\mu$
$$ F(a,z) = \int_{\R^d} \frac{d\mu(x)}{a\cdot x -z}, \ \ a \in R^n, \ z \in \C \setminus \R,$$
is developed. So far, we are not aware of the existence of such a
result.

As a direct application of the above discussion we can state the
following uniqueness criterion.

\begin{corollary} Let $\mu$ be a positive measure supported by the
positive octant $\R_+^d$. If all push-forward measures $\pi_\ell
(\mu)$ are determinate, where $\ell$ is an arbitrary semi-axis
contained in $\R_+^d$ passing through the origin and $\pi_\ell$
denotes the orthogonal projection onto $\ell$, then $\mu$ is
determinate.
\end{corollary}

The following section is devoted to an elaboration of the latter
observation.

\section{Disintegration techniques}

The present section focuses on disintegration of measures
techniques as another efficient tool for studying the determinacy
problem.

The starting point is the following {\it disintegration theorem}
for measures. Recall that $p(\nu) = p_* (\nu)$ denotes the
push-forward of a measure $\nu$ by a map $p$.

\begin{proposition}\label{disint} Suppose that $X$ and $T$ are a closed subsets of Euclidean
spaces and  $\nu$ is a finite positive Borel measure on $X$. Let
$p:X \to T$ be a $\nu$-measurable  mapping and let $\mu:=p(\nu)$.
Then there exist a mapping $t \to \lambda_t$ of $T$ into the set
of positive Borel measures on $X$ satisfying the following
conditions: \begin{itemize}
\item{\em (i)}~ supp~$\lambda_t \subseteq p^{-1}(t)$,\\
\item{\em (ii)}~ $\lambda_t(p^{-1}(t)) =1~~ \mu-{\rm a.e.}$,\\
\item{\em (iii)} $\int_X~f(x) ~d\nu(x) = \int_T ~d\mu(t) ~\int_X
~f(x)~d\lambda_t(x).$
\end{itemize}
\end{proposition}
This is a special case of Proposition 2.7.13 from \cite{Bou},
Chapter IX. Note that the measure $\nu$ is moderate, since $X$ and
$T$ are closed subsets of Euclidean spaces and hence locally
compact, and  the map $p$ is  $\nu$-proper, since $\nu$ is finite.
Let us retain the assumptions and notations of Proposition
\ref{disint}.

Let $\cA$ be a countably generated unital $\ast$-algebra of
$\nu$-integrable functions on $X$. Let $f \in \cA$ and $t \in T$.
We define linear functionals $L_t$ and $L$ on $\cA$ by
$$
L_t(f):= \int_X f(x)~d\lambda_t(x),~~L(f)\equiv L_\nu(f):= \int_X
f(x) ~d\nu(x).
$$
\begin{lemma}\label{l2ltf}
For $f\in \cA$, the function $t\to L_t(f)$ is in $L^2(T,\mu).$
\end{lemma}

Proof. We use freely the properties (i)--(iii) from the preceding
proposition. In particular, (iii) implies that $L_t(1)=1$. In view
of the above definitions and the Cauchy-Schwarz inequality we
obtain
\begin{align*}
&\int_T~ d\mu(t)~|L_t(f)|^2 = \int_T~ d\mu(t) ~|\int_X~ f ~d\lambda_t(x)|^2\leq\\
& \int_T ~d\mu(t) ( \int ~|f|^2~ d\lambda_t(x))~ ( \int~ 1~ d\lambda_t(x))=\\
& \int_T ~d\mu(t)~L_t(f\bar{f})L_t(1) = \int_T~d\mu(t)~
L_t(f\bar{f}) = L(f\bar{f}).~~~~~ \Box
\end{align*}

Suppose now that $\cB$ is a countably generated unital
$\ast$-algebra of functions on $T$ such that $ p^\ast \cB$ is a
$\ast$-subalgebra of $\cA$. That is all functions $(p^\ast f)(x)
:= f (p(x))$, where $f \in \cB$, belong to $\cA$ and the
involution commutes with the pull-back operation. Here and in what
follows we put $t{=}p(x)$ and write $\tilde{f}(t)=f(p(x))$. (Note
that $p^\ast \cB$ might be equal to $\cA$). We assume that the
following three conditions are
satisfied:\\

{\it (1) The measure $\mu$ is $\cB$-determinate (that is, if
$\mu^\prime$ is another measure on $T$ such that $\int ~d\mu(t)
{f}(t) = \int ~d\mu^\prime(t){f}(t)$  for all ${f}\in \cB$, then
$\mu=\mu^\prime$).

(2) $\cB$ is dense in $L^2(T,\mu)$.

(3) For $\mu$-almost all $t \in T$, the measure $\lambda_t$ is $\cA$-determinate on $p^{-1}(t)$ (that is, if $\lambda_t^\prime$ is another measure on the fibre $p^{-1}(t)$ such that $\int~f~d\lambda_t =\int~f~d\lambda_t^\prime $ for all $f \in \cA$, then $\lambda_t =\lambda_t^\prime$.) }\\

Our main result in this section is the following theorem.

\begin{thm}\label{genth}
Under the assumptions (1)-- (3) the measure $\nu$ is
$\cA$-determinate on $X$.
\end{thm}

Proof: Suppose that $\nu^\prime$ is another measure on $X$  such
that $\int~f~d\nu = \int ~f~d\nu^\prime$ for all $f \in \cA$.

Let $\mu{=}p(\nu)$,~ $\lambda_t$ and $\mu\prime{=}p(\nu^\prime)$,~
$\lambda_t^\prime$, respectively, be the corresponding measures
from the disintegration theorem and let $L_t$ and
$L_t^\prime$,respectively, be the corresponding linear functionals
for $\nu$ and $\nu^\prime$, respectively. For $\tilde{f} \in \cB$
we compute
\begin{align*}
&\int ~d\mu^\prime(t)~ \tilde{f}(t)=
\int~d\mu^\prime(t)\tilde{f}(t) \int ~d\lambda_t^\prime(x) =
 \int~ d\mu^\prime(t) \int f(p(x))~d\lambda_t^\prime(x) =\\
&\int f(p(x)) d\nu^\prime(x) =L_{\nu^\prime}(f(p(x)))=
L_\nu(f(p(x))) = \int ~d\mu(t) \tilde{f}(t).
\end{align*}
Since $\mu$ is $\cB$-determinate by assumption (1), it follows
that $\mu=\mu^\prime$.

Since $\cB$ is countably generated and dense in $L^2(T,\mu)$ by
assumption (2), there are functions $\varphi_n =
p^\ast(\tilde{\varphi_n}) \in \cA$, $n {\in} N$, such that
$\{\tilde{\varphi}_n;n {\in} N\}$ is an orthonormal basis of
$L^2(T,\mu)$. Fix $f \in \cA$. We compute the Fourier coefficents
of the function $L_t(f)$ of $L^2(T,\mu)$ with respect to this
orthonormal basis by
\begin{align*}
&\int_T ~d\mu(t) \tilde{\varphi}_n(t)L_t(f) = \int_T
d\mu(t)~\tilde{\varphi}_n(t) \int f(x)~d\lambda_t(x)
 =\\
& \int~d\mu(t)~\int \varphi(p(x))f(x)~d\lambda_t(x)
=L(\varphi_n(p(x))f(x)).
\end{align*}
Since $\mu=\mu^\prime$, the same reasoning with $\lambda_t$
replaced by $\lambda_t^\prime$ shows that
\begin{align*}
\int_T ~d\mu(t)\tilde{\varphi}_n(t) L_t^\prime(f) =
L(\varphi_n(p(x))f(x)).
\end{align*}
Therefore, both functions $L_t(f)$ and $L_t^\prime(f)$ from
$L^2(T,\mu)$ (by  Lemma \ref{l2ltf}) have the same developments
\begin{align*}
L_t(f)=\sum_n~ L( \varphi_n(p)f) \tilde{\varphi}_n(t),~~
L_t^\prime(f)= \sum_n~ L(\varphi(p)f) \tilde{\varphi}_n(t)
\end{align*}
in $L^2(T,\mu)$. Consequently, we have $L_t(f) =L_t^\prime(f)$
$\mu$-a-e. on $T$. That is, there is a $\mu$-null subset $M_f$ of
$T$ such that $\int~ f(x)~d\lambda_t(x)=\int
f(x)~d\lambda_t^\prime(x)$ for  $t \in T\backslash M_f$. Since
$\cA$ is countably generated, there is a $\mu$-null subset $M$ of
$T$ such that the latter holds for {\it all} $f \in \cA$ and for
$t \in T\backslash M$. From assumption (3) it follows therefore
that $\lambda_t=\lambda_t^\prime$ $\mu$-a.e. on $T$. From the
disintegration formula (iii) we obtain $\nu = \int d\mu(t)~
\lambda_t = \int d\mu^\prime~ \lambda_t^\prime =\nu^\prime$.
\hfill $\Box$\\

\noindent We now specialize the preceding general theorem to the
moment problem and derive  a "reduction procedure" for proving
determinacy.

Let $X$ be a closed subset of $\R^d$, $\cA:=\C[x_1,\dots,x_d]$ and
$\cB:=\C[t_1,\dots,,t_m]$. Let $p_1,\dots,p_m \in
\R[x_1,\dots,x_d]$ be  polynomials and define a mapping $p:X\to T$
by $p(x)=(p_1(x),\dots,p_m(x))$, where $T$ is a closed subset of
$\R^m$ such that $p(X)\subseteq T$. Note that in the case $X=\R^m$
the fibres are just the real algebraic varieties
$$ p^{-1}(t)=\{x \in \R^n: p_1(x)=t_1,\dots,p_k(x)=t_k \},~~~ t\in T.
$$
Suppose that $\nu \in \cM(X)$. Let $\mu$ and $\lambda_t$ denote
the corresponding measures from Proposition \ref{disint}. Since
$L_\nu(f(p(x)) =\int ~d\mu(t) \tilde{f}(t) < \infty$ for all
$\tilde{f} \in \R[t_1,\dots,t_m]$, we also have that $\mu \in
\cM(X)$. From Lemma \ref{l2ltf} it follows that $\lambda_t \in
\cM(X)$ $\mu$-a. e. on $T$. Note that each measure $\lambda_t$ is
supported on the fibre $p^{-1}(t)$. Recall that a measure $\mu \in
\cM(T)$ is {\it strictly determinate} if it is determinate on $T$
and if the polynomials $C[t_1,\dots,t_m]$ are dense in $L^2(\mu)$.

In the preceding setup Theorem \ref{genth} can be restated as
follows.

\begin{thm}\label{applth}
If $\mu$ is strictly determinate  on $T$ and $\lambda_t$ is
determinate on the fibre $p^{-1}(t)$ for $\mu$-almost all $t \in
T$, then $\nu$ is determinate on $X$.
\end{thm}

Theorem \ref{applth} combined with the fact that measures with
compact support are always (strictly) determinate gives the
following two corollaries.

\begin{corollary}\label{corbound1}
If $\mu$ is strictly determinate on $T$ and the fibre $p^{-1}(t)$
is bounded for $\mu$-almost all $t \in T$, then $\nu$ is
determinate on $X$.
\end{corollary}

\begin{corollary}\label{corbound2}
If $T$ is compact and the measure $\lambda_t$ is determinate on
the fibre $p^{-1}(t)$ for $\mu$-almost all $t \in T$, then $\nu$
is determinate on $X$.
\end{corollary}

There is a large number of applications of  Theorem \ref{applth}
and Corollaries \ref{corbound1} and \ref{corbound2} by specifying
the set $X$ and the polynomials $p_1,\dots,p_m$ occuring therein.
We mention three such results and retain the notations and the
setup introduced before Theorem \ref{applth}.\\

{\bf 1.} Let $p_1(x)=x_1,\dots,p_m(x)=x_m$, $m<d$, and let $X$ and
$T$ be closed subsets of $\R^d$ and $\R^m$, respectively, such
that $p(X)\subseteq T$. Then Theorem \ref{applth} yields:\\

 {\it The measure $\nu$ is determinate on $X$ if $\mu=p(\nu)$ is strictly
 determinate on $T$ for $\C[t_1,\cdots,t_k]$ and the fibre measure $\lambda_{(t_{1},\cdots,t_m)}$
 is determinate on $p^{-1}(t)$ for $\C[x_1,\cdots,x_d]$ and $\mu$-almost all $t\in
 T$.}\\

\noindent {\bf 2.} Set $p(x)=x_1^2+\cdots+x_n^2$. Let $X=\R^n$ and
$T:=p(X)=[0,\infty)$. Since strict determinacy is the same as
determinacy in dimension one, Corollary \ref{corbound1} gives:\\

 {\it The measure $\nu$ is determinate
if $\mu=p(\nu)$ is determinate on $[0,\infty)$}.\\

\noindent For rotation invariant measures on $\R^d$ this assertion
has been obtained in \cite{Bergth}.\\

\noindent {\bf 3.} Suppose that
$p_1,\dots,p_m\in\R[x_1,\dots,x_d]$ are polynomials which are {\it
bounded} on the closed subset $X$ of $\R^d$. Put
$$
\alpha_j={\rm inf}~\{p_j(x); x \in X\}, \beta_j={\rm sup}~
\{p_j(x); x \in X \}, T=[\alpha_1,\beta_1]\times \cdots \times
[\alpha_m,\beta_m].
$$
An immediate consequence of Corollary \ref{corbound2} is the
following assertion:\\

{\it The measure $\nu$ is determinate on $X$ if the fibre measures
$\lambda_t$ are determinate on $p^{-1}(t)$ for $\mu$-almost all
$t\in T$.}\\

\noindent If the closed set $X$  has non-constant bounded
polynomials, it is often possible to reduce question about the
moment problem on $X$ to their counter--parts for the moment
problem on their (in general lower dimensional) fibre sets. An
existence theorem of this kind was proved in \cite{S2}. The
preceding assertion is a similar result for the determinacy
question.

\section{Geometric determinateness}
We reproduce below in a condensed and almost identical form a
result from \cite{PS} which illustrates the thesis that the
geometry of the support of a measure implies, under specific
conditions, its determinateness. The notation and terminology,
standard in real algebraic geometry, are explained in the same
article \cite{PS} or in \cite{Sch}.

The simplest example is of course a compact support, in which case
Stone-Weierstrass theorem leads immediately to the uniqueness
conclusion. The next proposition generalizes this observation.

\begin{proposition}\label{Hseparatespts}
Let $X$ be an affine $\R$-variety, and let $K$ be a closed subset
of $X(\R)$. If the algebra ${\mathcal H}(K)=\{p\in\R[X]\colon p$
is bounded on $K\}$ separates the points of $K$, then every
$K$-moment problem is determinate.
\end{proposition}

\begin{proof}
First assume that $H={\mathcal H}(K)$ is generated by finitely
many elements $h_1,\dots,h_m$ as an $\R$-algebra. The map
$h:=(h_1,\dots,h_m)\colon K\to\R^m$ is injective, and the subset
$\ol{h(K)}$ of $\R^m$ is compact. From Theorem \ref{pedcor1} we
infer that every $K$-moment problem is determinate. The situation
when ${\mathcal H}(K)$ is not finitely generated can be reduced to
the previous case via a finite chain condition and an inductive
limit argument, see for details \cite{PS}.\end{proof}

We are now discussing cases to which Proposition
\ref{Hseparatespts} applies by first considering one dimensional
varieties.

Let $X$ be an affine curve over $\R$, and let $\ol X$ be its
(good) completion. That is, the unique (up to unique isomorphism)
projective curve which contains $X$ as a Zariski dense open subset
and whose points in the complement of $X$ are nonsingular. Let
$S=\ol X-X$ (a finite set), and let $K$ be a closed subset of
$X(\R)$. For simplicity, assume that $X$ is irreducible. Following
\cite{Sch} we'll say that $K$ is \emph{virtually compact} if $S$
contains at least one point which is either non-real or does not
lie in the closure $\ol K$, the closure being taken in $\ol
X(\R)$.

Let $\mathcal H={\mathcal H}(K)$ be the subring of $\R[X]$
consisting of all regular functions which are bounded on $K$.
Regarding elements $p\in\R[X]$ as rational functions on $\ol X$,
$p$ is bounded on $K$ if and only if none of the points of $\ol
K\cap S$ is a pole of $p$. So $\mathcal H=\mathcal{O} (\ol X-T)$,
where $T$ is the set of points in $S$ which do not lie in $\ol K$.
Therefore we see (cf.\ \cite{Sch}, Lemma 5.3) that $K$ is
virtually compact if and only if $\mathcal H\ne\R$, and that in
this case $\mathcal H$ separates the points of $X(\R)$. Hence:

\begin{thm}\label{virtcptmpdet}
Let $X$ be an irreducible affine curve over $\R$, and let $K$ be a
closed subset of $X(\R)$. If $K$ is virtually compact then every
moment problem on $K$ is determinate. \qed
\end{thm}

The condition that $X$ is irreducible can be removed (see
\cite{Sch}, Definition 5.1 and Lemma 5.3).

For the case of one-dimensional sets $K$, and for the
determinateness question, this leaves us with the case where $K$
is not virtually compact. In other words, the case where every
polynomial which is bounded on $K$ is constant on $K$. Going back
to Proposition \ref{suffdetermined}, the following observation is
in order.

\begin{lemma}
Assume that the affine curve $X$ is irreducible, and let $K\subset
X(\R)$ be a closed subset. Then $\mathcal A_0(K,f)$ separates the
points of $X(\R)$ for every non-constant $f$ in $\R[X]$ with
$f\ge1$ on $X(\R)$.
\end{lemma}

(Instead of $f\ge1$ it is only needed here that $f$ vanishes
nowhere on $X(\R)$.)

\begin{proof}
Let $Y\subset\ol X$ be the open set where $f$ is regular. So $Y$
is affine, contains $X$, and the points in $\ol X-Y$ are
nonsingular on $\ol X$. Since the rational function $\frac1f$
vanishes in the points of $\ol X-Y$, there exists for every
$q\in\R[Y]$ an integer $k\ge1$ such that $\frac q{f^k}$ vanishes
at all points of $\ol X-Y$. Let
$$I=\bigl\{q\in\R[Y]\colon\forall y\in Y-X,\ q(y)=0\bigr\},$$
an ideal of $\R[Y]$. If $q\in I$, and if $k\ge1$ is chosen for $q$
as before, the rational function $\frac q{f^k}$ lies in $\mathcal
A_0(K,f)$. Since the elements of $I$ separate the points of
$X(\R)$, the lemma follows.
\end{proof}

Consequently the following application is obtained.

\begin{corollary}
Let $X$ be an irreducible affine curve and $K\subset X(\R)$ a
closed set. If there are a non-constant $f\in\R[X]$ with $f\ge1$
on $K$ and a sequence $p_n$ in $\R[X]$ with $fp_n\to1$ under
$||\cdot||_{L,2}$, the $K$-moment problem $L$ is determinate. \qed
\end{corollary}

There are higher-dimensional cases as well which are non-compact
and to which Proposition \ref{Hseparatespts} applies. Here is a
class of examples.

\begin{example}
Let $K_1$ be a compact subset of $\R^n$, let $f\colon K_1\to\R$ be
a continuous function, and let $K=\{(x,t)\in K_1\times\R\colon
tf(x)= 1\}$, a closed subset of $\R^{n+1}$. Then ${\mathcal H}(K)$
separates the points of $K$. By Proposition \ref{Hseparatespts},
therefore, any $K$-moment problem is determinate.
\end{example}

To illustrate Theorem \ref{virtcptmpdet}, we give a couple of
examples of one-dimensional sets which are virtually compact but
not compact. For simplicity we stick to subsets of the plane. Let
$p\in \R[x,y]$ be an irreducible polynomial, let $X$ denote the
plane affine curve $p=0$.\\

\emph{If the leading form (i.e., highest degree form) of $p(x,y)$
is not a product of linear real factors, then every closed subset
$K$ of $X(\R)$ is virtually compact. Thus, every $K$-moment
problem is determinate.}\\

But also if the leading form of $p$ is a product of real linear
forms, $X(\R)$ may be virtually compact (let alone closed subsets
$K$). The reason is that, although the Zariski closure of the
curve $X$ in the projective plane contains only real points at
infinity, some of them may be singular and may blow up to one or
more non-real points. An example is given by the curve $p=
x+xy^2+y^4=0$.

Finally, even if the entire curve $X(\R)$ itself fails to be
virtually compact, suitable non-compact closed subsets $K$ may
still be. For example, this is so for the hyperelliptic curves
$y^2=q(x)$, where $q$ is monic of even degree, not a square. For
example, one easily checks that if $\deg(q)$ is divisible by~$4$,
then $K$ is virtually compact if (and only if) the coordinate
function $y$ is bounded on $K$ from above or from below.

\end{document}